\documentclass[10pt]{article}
\usepackage{amssymb}
\usepackage{amsmath}
\usepackage{amsthm}
\theoremstyle{plain}
\newtheorem{theorem}{Theorem}[section]
\newtheorem{proposition}[theorem]{Proposition}
\newtheorem{lemma}[theorem]{Lemma}
\newtheorem{corollary}[theorem]{Corollary}
\theoremstyle{definition}
\newtheorem{definition}[theorem]{Definition}
\newtheorem{remark}[theorem]{Remark}
\newtheorem{example}[theorem]{Example}

\theoremstyle{remark}

\newcommand{\N}{\mathbb{N}}
\newcommand{\I}{{\bf I}}
\newcommand{\Z}{\mathbb{Z}}

\newcommand{\C}{\mathbb{C}}

\newcommand{\T}{\mathbb{T}}

\newcommand{\Zj}{{\cal Z}}

\newcommand{\noteq}{=\hspace{-.60em}/}
\DeclareMathOperator{\id}{id}

\DeclareMathOperator{\Aut}{Aut}
\DeclareMathOperator{\Ad}{Ad}
\DeclareMathOperator{\WInn}{WInn}

\DeclareMathOperator{\Tr}{Tr}

\newcommand{\bprf}{\noindent{\it Proof.\ }}

\newcommand{\eprf}{\hspace*{\fill} \rule{1.6mm}{3.2mm} \vspace{1.6mm}}
\newcommand{\benu}{\begin{enumerate}\renewcommand{\labelenumi}{{\rm (\roman{enumi})}}\renewcommand{\itemsep}{0pt}}
\newcommand{\eenu}{\end{enumerate}}

\begin{document}
\title{ Certain aperiodic automorphisms of unital simple projectionless  $C^*$-algebras}
\author{ Yasuhiko Sato \\}
\date{\small  Department of Mathematics, Hokkaido University, \\ Sapporo 060-0810, Japan \\ e-mail : s053018@math.sci.hokudai.ac.jp}

\maketitle

\begin{abstract}
 Let $G$ be an inductive limit of finite cyclic groups and let $A$ be a unital simple projectionless $C^*$-algebra with $K_1(A) \cong G$ and with a unique tracial state, as constructed based on dimension drop algebras by Jiang and Su. First, we show that any two aperiodic elements in $
\Aut(A)/\WInn(A)$ are conjugate, where $\WInn(A)$ means the subgroup of $\Aut(A)$ consisting of automorphisms which are inner in the tracial representation.

In the second part of this paper, we consider a class of unital simple $C^*$-algebras with a unique tracial state which contains the class of unital simple $A\T$-algebras of real rank zero with a unique tracial state. This class is closed under inductive limits and under crossed products by actions of $\Z$ with the Rohlin property. Let $A$ be a TAF-algebra in this class. We show that for any automorphism $\alpha$ of $A$ there exists an automorphism $\widetilde{\alpha}$ of $A$ with the Rohlin property such that $\widetilde{\alpha}$ and $\alpha$ are asymptotically unitarily equivalent. In its proof we use an aperiodic automorphism of the Jiang-Su algebra.   


\end{abstract}
\section{Introduction}\label{INTRO}
The classification theory of automorphisms of von Neumann algebras was studied by Connes \cite{Co} using the Rohlin property for automorphisms of von Neumann algebras. To be specific, he classified automorphisms of the injective type $I\hspace{-.1em}I_1$-factor up to outer conjugacy. In a case of finite $C^*$-algebras, Kishimoto showed some results (for example, \cite{Kis1}, \cite{Kis}, \cite{Kis3}, \cite{Kis2}) for classification of automorphisms with the Rohlin property adapted to $C^*$-algebras, and in a case of infinite $C^*$-algebras, Nakamura \cite{Na} classified aperiodic automorphisms of nuclear purely infinite simple $C^*$-algebras by KK-class proving the Rohlin property for such automorphisms. In order to define the Rohlin properties, these previous works considered only $C^*$-algebras of real rank zero. In this paper, we study aperiodic automorphisms of certain unital simple projectionless $C^*$-algebras with a unique tracial state, which were constructed by Jiang and Su \cite{JS}. Instead of the Rohlin property we require the automorphisms to satisfy that any non-zero power is not weakly inner in the unique tracial representation. Our first main result Theorem \ref{thm:wa} says that: if two automorphisms of such a unital simple projectionless $C^*$-algebra $A$ satisfy the above property, they are conjugate in the quotient group $\Aut(A)/ \WInn (A)$, where 
\[ \WInn (A) = \{ \alpha \in \Aut (A) :\pi_{\tau}\circ \alpha= \Ad W \circ \pi_{\tau} ,\quad W\in U(\pi_{\tau}(A)'')\}.\]
Here $\tau$ is the unique tracial state of $A$ and $\pi_{\tau}$ is the GNS representation associated with $\tau$. Note that we do not require the two automorphisms to have the same KK-class and that an automorphism which is weakly inner in the tracial representation need not have the same KK-class as the identity automorphism.

In the second half, we consider an existence problem  of automorphisms with the Rohlin property for unital simple TAF-algebras, where the class of TAF-algebras is introduced by H.Lin \cite{Lin2} and is shown to include the simple unital AH-algebras of real rank zero with slow dimension growth (Elliott and Gong \cite{EG}). We introduce a technical property called (SI) (where SI stands for  small isometry for a $C^*$-algebra) which roughly says that if two central sequences of projections are given such that one of them is infinitesimally small compared to the other in the sequence algebra, then in fact so in the central sequence algebra, and we mainly consider a class of unital simple TAF-algebras with the property (SI) and with a unique tracial state. The property (SI) is devised to satisfy that the class of unital simple $C^*$-algebras with (SI) includes the unital simple A$\T$-algebras of real rank zero with a unique tracial state and is closed under taking crossed products by automorphisms with the Rohlin property and taking inductive limits. Our second result Theorem \ref{thm:1} says that if a unital TAF-algebra $A$ absorbs the Jiang-Su algebra, and has the property (SI), then any automorphism of $A$ is asymptotically unitarily equivalent to a Rohlin automorphism.

The contents of this paper are as follows. In Section 2, we reconstruct a certain projectionless $C^*$-algebras in Lemma \ref{lem:igh} by an adaptation of the construction of the Jiang-Su algebra, and we prove the first main result Theorem \ref{thm:wa} by using this construction and an intertwining argument. In Section 3, we introduce the definition of the property (SI), and we show in Proposition \ref{Prop:AT} that any unital simple A$\T$-algebra  of real rank zero satisfies the property (SI). We also show the permanence properties of the class Proposition \ref{prop:crossedproduct} and Proposition \ref{prop:inductivelimit}. In Section 4, we show the second main results Theorem \ref{thm:1}.

Concluding this section, we prepare some notations. When $A$ is a $C^*$-algebra, we denote by $P(A)$ the set of projections in $A$, $U(A)$ the unitary group of $A$, $A_+$ the set of positive elements in $A$, $A^1$ the set $\{a\in A: \|a\| \leq 1 \}$, $T(A)$ the  tracial state space of $A$, $\Aut (A)$ the automorphism group of $A$, $A^{\infty}$ the quotient $\ell^{\infty}(\N,A)/c_0(A)$, and $A_{\infty}$
 the relative commutant $A^{\infty} \cap A'.$
We define an inner automorphism of $A$ by $\Ad u(a)=uau^*$ for $u \in U(A)$ and $a\in A$. We use $M_n$ to denote the $C^*$-algebra of $n\times n$ matrices with complex entries.
 A dimension drop algebra means  a 
$C^*$-algebra of the form $\I [m_0,m,m_1]=$
\[\{f \in C([0,1],M_m); f(0) \in M_{m_0}\otimes 1_{m/m_0},\ f(1)\in 1_{m
/m_1}\otimes M_{m_1}\},\]
where $m_0$ and $m_1$ are divisors of $m$ (see Definition 2.1 in \cite{JS}). We denote by $(m,n)$  the greatest common divisor of $m$, $n$.
When $A$ is a $C^*$-algebra and $B$ is a finite dimensional $C^*$-subalgebra of $A$ such that $B= \bigoplus_{n=1}^{N} M_{d_n}$,  we define a conditional expectation $\Phi_{B} : A \rightarrow A\cap B'$ by 
\[\Phi_{B}(a)=\sum_{n=1}^{N} d_{n}^{-1}\sum_{1\leq i,j \leq d_n} e_{i,j}^{(n)}ae_{j,i}^{(n)},\quad a \in A,\]
where $\{e^{(n)}_{i,j}\}$ is a system of matrix units of $B$. 
\section{The weakly aperiodic automorphism of projectionless $C^*$-algebras}\label{sec:W}
In this section, first we reconstruct a unital simple projectionless $C^*$-algebra $A$ such that $A$ has a unique tracial state $\tau$ and $K_1(A)=G$, where $G$ is any abelian group obtained as the inductive limit of a sequence of finite cyclic groups. We employ a different method from the one used by Jiang and Su in the proof of Theorem 4.5 \cite{JS} so that we obtain a UHF algebra $B$ such that $\pi_{\tau}(A) \subset B \subset \pi_{\tau}(A)''$. Note that, by the uniqueness result in \cite{JS}, the new $A$ is isomorphic to the original one in \cite{JS}.  Using this $B$ with some additional properties we then obtain the weakly outer conjugacy result in Theorem \ref{thm:wa}.

For $ m, n, r_0, r_1  \in \N$ with $m>r_0+r_1 $ and $s\in [0,1]$, we define continuous paths $\xi_{j,s}$ in $[0,1]$ (see \cite{JS}) by
\begin{eqnarray}
\xi_{j,s}(t)=\left\{ \begin{array}{ll}
st,\quad & 1\leq j \leq r_0 \\
s,\quad  & r_0< j \leq  m-r_1 \\
(1-s)t+s,\quad & m-r_1<j\leq m, \\
\end{array} \right.\nonumber
\end{eqnarray}  
with $t\in [0,1]$.  
And we define an injective $*$-homomorphism $\xi_{n,m,r_0,r_1,s}:C([0,1])
\otimes M_n \rightarrow C([0,1])\otimes M_{mn}$ as 
\[\xi_{n,m,r_0,r_1,s}(f)=f\circ\xi_{1,s} \oplus f\circ\xi_{2,s} \oplus \cdots
 \oplus f\circ\xi_{m,s},\quad f \in C([0,1])\otimes M_n.\]

The proof of the following lemma is an adaptation of the proof of Proposition 2.5 \cite{JS}.
\begin{lemma}\label{lem:igh}
Let $G_n$ be a finite cyclic group $\Z_{g_n} $ for each $n\in\N$ and let $\gamma_n$ be an 
injective homomorphism from $G_n$ into $G_{n+1}$ such that $\gamma_n(1_n)
=c_{n+1}1_{n+1}$, where $1_n=[1] \in \Z_{g_n} $ and $g_n$, $g_{n+1}$, and $c_{n+1}$ are natural numbers such that $c_{n+1}g_n=g_{n+1} $. 
For any natural numbers $p_n, q_n,k$ such that $p_n $ and $q_n $ are relatively
prime and $ k >g_n(p_n+q_n) $, there exist natural numbers $p_{n+1}, q_{n+1}$, a unitary $u_n \in U(C([0,1])\otimes M_{d_{n+1}}) $, and
 a unital injective $*$-homomorphism $\varphi_n$ from $\I[p_n, d_n,q_n] $ into 
$\I [p_{n+1},d_{n+1}, q_{n+1}] $ such that
\[(p_{n+1}, q_{n+1})=1,\quad p_n|p_{n+1},\quad q_n|q_{n+1},\quad p_{n+1}/p_n>k,\quad  q_{n+1}/q_n>k, \]
\[\varphi_n = \Ad u_n \circ (\xi_{d_n,m_n,r_{n,0},r_{n,1},1/2}\otimes 1_{M_{c_{n+1}}}),
\quad (\varphi_n)_*=\gamma_n, \]
where $d_n=p_ng_nq_n$, $m_n=p_{n+1}q_{n+1}/(p_nq_n)$, $r_{n,0}$ and $r_{n,1}$ are some natural numbers, and $(\varphi_n)_*$ is the induced map on $K_1$-group.  
\end{lemma}
\bprf
We obtain a natural number $k_0$ such that $(k_0,q_n)=1 $, $ k_0\equiv1
\pmod{g_n}$, and $k_0>k$ (e.g., $k_0=ag_nq_n+1,\ a\in \N$), and we define $p_{n+1}$ as $p_nk_0$.
Similarly we obtain a natural number $k_1$ such that $(k_1,p_{n+1})=1$, 
$k_1 \equiv 1 \pmod{g_n} $, and $k_1>k_0$ and we define $q_{n+1}$ as $q_nk_1$.
Let $r_{n,0}$, $r_{n,1}$ be natural numbers such that 
\begin{eqnarray}
k_0k_1&=&ag_nq_{n+1}+r_{n,0} \nonumber \\
&=&bg_np_{n+1}+r_{n,1}, \nonumber
\end{eqnarray}  
where $a,b \in \N $, $0< r_{n,0} \leq g_nq_{n+1} $, and $0<r_{n,1} \leq g_np_{n+1}$.
Remark that $r_{n,0} \equiv r_{n,1} \equiv 1 \pmod{g_n}$ and $(p_{n+1},q_{n+1})=1$. Put $m_n=k_0k_1$. By $k_1>k_0>k>g_n(p_n+q_n)$, we have that
\[m_n-(r_{n,0}+r_{n,1})>k_1g_n(p_n+q_n)-(g_nq_nk_1+g_np_nk_0)=g_np_n(k_1-k_0)>0.\]

We denote $\xi_{d_n,m_n,r_{n,0},r_{n,1},s}$ by $\xi_s $. Since 
\[ \displaystyle \xi_s(f)(0)=\bigoplus_{i=1}^{r_{n,0}} f(0) \oplus \bigoplus_{i=1+r_{n,0}}^{m_n}f(s)\] and $f(0) \in M_{p_n} \otimes 1_{g_nq_n} $ for any $f \in \I[p_n,d_n,q_n] $, $g_nq_{n+1}|g_nq_nr_{n,0}$, and $ g_nq_{n+1}|$ $m_n-r_{n,0} $, we obtain $u_0 \in 
U(M_{m_nd_n})$ independent of $s$ such that 
\[u_0 \cdot \xi_s(f)(0) \cdot u_0^* \in M_{p_{n+1}} \otimes 1_{g_nq_{n+1}},\]
for any $s \in [0,1]$ and $f \in \I [p_n,d_n,q_n] $. 
 Similarly, since 
\[\displaystyle \xi_s(f)(1)=\bigoplus_{i=1}^{m_n-r_{n,1}} f(s) \oplus \bigoplus_{i=1+m_n-r_{n,1}}^{m_n}f(1)\] and $f(1) \in 1_{g_np_n} \otimes M_{q_n} $ for any $f \in \I[p_n,d_n,q_n] $, $g_np_{n+1}|g_np_nr_{n,1}$, and $ g_np_{n+1}|$ $m_n-r_{n,1} $, we obtain $u_1 \in U(M_{m_nd_n}) $ independent of $s$ such that 
\[u_1 \cdot \xi_s(f)(1) \cdot u_1^* \in 
1_{g_np_{n+1}} \otimes M_{q_{n+1}},\]
for any $s \in [0,1]$ and $f \in \I [p_n,d_n,q_n] $. Then we obtain a unitary $u \in U(C([0,1])\otimes M_{m_nd_n})$ such that $u(0)=u_0$ and $u(1)=u_1$. 

We define an injective unital $*$-homomorphism $\varphi_n':C([0,1])\otimes M_{d_n} \rightarrow 
C([0,1])\otimes M_{m_nd_n} $ by 
\[\varphi_n'(f)= \Ad u \circ \xi_{1/2}(f),\quad f \in C([0,1])
\otimes M_{d_n}. \]
Then $\varphi_n'$ satisfies $\varphi_n'
(\I[p_n,d_n,q_n]) \subset \I[p_{n+1},m_nd_n,q_{n+1}]$ by the definition of unitaries $u_0, u_1 $.
Let $f_1$ be a unitary in $\I[p_n,d_n,q_n] $ such that \[f_1(t)=\exp(2\pi it)e_n
 + (1_{d_n}-e_n), \] where $e_n$ is a minimal projection of $M_{d_n}$. We identify $[f_1]\in K_1(\I[p_n,d_n,q_n])$ as a generator $1_n \in \Z_{g_n}$. Since $u$ is independent of $s$, it follows that $\varphi_n'(f_1)=u\xi_{1/2}(f_1)u^*$ and $u\xi_{0}(f_1) u^* $ are homotopic in $U(\I [p_{n+1},m_nd_n,q_{n+1}]) $. By $\xi_{0}(f_1)(0)=1_{m_nd_n}=\xi_{0}(f_1)(1)$, $u\xi_{0}(f_1) u^* $ and $\xi_{0}(f_1)$ are homotopic in $U(\I[p_{n+1},m_nd_n,q_{n+1}]). $
Thus we can conclude that  
\[ [\varphi_n'(f_1)]_1=[\xi_{0}(f_1)]_1 = r_{n,1} 1_n = 1_n. \]

We define a unitary $u_n$ in $C([0,1])\otimes M_{d_{n+1}} $ as $u\otimes 1_{c_{n+1}}$ and define an injective unital $*$-homomorphism $\varphi_n:C([0,1])\otimes M_{d_n} \rightarrow C([0,1])\otimes M_{d_{n+1}}$ by 
\[\varphi_n= \Ad u_n \circ (\xi_{1/2}\otimes 1_{c_{n+1}}). \] 
This $\varphi_n$ satisfies the required conditions 
$\varphi_n(\I[p_n,d_n,q_n]) \subset \I[p_{n+1},d_{n+1},q_{n+1}]$ and $(\varphi_n)_*(1_n) = c_{n+1}1_{n+1}$ in $ K_1(\I[p_{n+1},d_{n+1},q_{n+1}])$.\eprf

\begin{lemma}\label{lem:ifc}
Let $G$ be an abelian group which is the inductive limit of a sequence of finite cyclic groups. Then there exists a sequence of dimension drop algebras $A_n=\I[p_n,d_n,q_n]$ with $(p_n,q_n)=1$ and injective $*$-homomorphisms $\varphi_n:A_n\rightarrow A_{n+1}$ satisfying the following conditions: the inductive limit $C^*$-algebra $\displaystyle A=\lim_{\longrightarrow}(A_n,\varphi_n)$ is a unital simple projectionless $C^*$-algebra with a unique tracial state $\tau $ such that $K_1(A)=G$ and the probability measure on $[0,1]$ determined by $\tau|_{A_n}$ is faithful and has support in $(0,1)$ for any $n \in \N$, where $A_n$ is canonically regarded as a subalgebra of $A$.     
\end{lemma}
\bprf Let $G_n$ be the finite cyclic group $\Z_{g_n}$ of order $g_n$ and let $\gamma_n:G_n \rightarrow G_{n+1} $ be a connecting morphism such that $\displaystyle G=\lim_{\longrightarrow}(G_n,\gamma_n)$. We may assume that $\gamma_n $ is an injective map. For $\gamma_n$ there exists a natural number $x$ such that $\gamma_n(1_n)=xc_{n+1}1_{n+1}$, where $g_{n+1}=c_{n+1}g_n$ and $1_n=[1]\in G_n$. Since $\gamma_n$ being injective implies that $(x,g_n)=1$, we obtain that $\gamma_x:G_n \rightarrow G_n$ with $\gamma_x(1_n)=x1_n$ is an isomorphism of $G_n$. Let $\gamma_n':G_n \rightarrow G_{n+1}$ be a group homomorphism such that $\gamma_n'(1_n)=c_{n+1}1_{n+1}$. Since $\gamma_n=\gamma_n'\circ \gamma_x$, $G$ is isomorphic to the inductive limit $\displaystyle \lim_{\longrightarrow}(G_n,\gamma_n')$ and thus we may assume that $\gamma_n(1_n)=c_{n+1}1_{n+1} $. 

By the above Lemma \ref{lem:igh} we can inductively obtain dimension drop algebras $A_n$ and $*$-homomorphisms $\varphi_n:A_n
\rightarrow A_{n+1}$ 
such that $K_1(A_n)=G_n$, $(\varphi_n)_*=\gamma_n $. Thus the inductive limit $C^*$-algebra $\displaystyle A=\lim_{\longrightarrow}(A_n,\varphi_n)$ satisfies that $K_1(A)=G$. By the construction of $\varphi_n$ and Lemma 2.7 of \cite{JS}, $A$ is a simple $C^*$-algebra. By the proof of the Proposition 2.8 of \cite{JS}, $A$ has a unique tracial state $\tau$.

Define a probability measure $\mu_n$ on $[0,1]$ by 
\[ \tau \widetilde{\varphi}_n (f) = \int_{[0,1]} \Tr_{d_n} (f(t)) d \mu_n (t),\]for $f \in A_n$, where $\widetilde{\varphi}_n $ is the canonical embedding of $A_n$ into $A$. We have to show that $\mu_n (\{ 0\})= 0 $ and $\mu_n (\{ 1\})=0$.

Assume that $\mu_n (\{0\}) >0$. Since
\[ \tau\widetilde{\varphi}_n(f)= \tau \widetilde{\varphi}_{n+1} \varphi_n (f) =
\int _{[0,1]} \Tr_{d_{n+1}} (\xi^{(n)} (f) (t)) d \mu_{n+1} (t),\]
where $\xi^{(n)} = \xi_{1/2} \otimes 1_{c_{n+1}}$ in the proof of Lemma \ref{lem:igh}, and since $f_n(0)$ appears in $\xi^{(n)} (f) (t)$ only at $t=0$ as $c_{n+1}r_{n,0}$ copies among the $c_{n+1}m_n$ direct summands, it follows that 
\[ \mu_{n+1} (\{ 0\}) \frac{r_{n,0}}{m_n} = \mu_n ( \{ 0\}).\] 
Since $\frac{r_{n,0}}{m_n} \leq \frac{1}{2}$, this would lead to $\mu_{n+m}(\{0\}) >1 $ for some $m$, a contradiction. Hence it follows that $\mu_n(\{0\}) =0 $, and similarly, that $\mu_n(\{1\})=0$.
  \eprf

The next proposition is based on the Thomsen's theorem \cite{Thom}.
\begin{proposition}\label{prop:eUHF}
Let $G$ be the inductive limit of a sequence of finite cyclic groups. Let $A$ be a simple projectionless $C^*$-algebra which is obtained in Lemma \ref{lem:ifc} for $G$, $\tau$ the unique tracial state of $A$, and $\pi_{\tau}$ the GNS representation associated with $\tau$. Then there exists a UHF $C^*$-subalgebra $B$ of $\pi_{\tau} (A)'' $ and an increasing sequence $(B_n)$ of matrix $C^*$-subalgebras of $B$ such that $\pi_{\tau} (A) \subset B$, $1_{B_n}=1_{B}$, $\displaystyle B= \overline{(\bigcup B_n)}^{\|\cdot\|} $, and 
\[(\pi_{\tau}(A) \cap B_n')''= \pi_{\tau}(A)'' \cap B_n'.\]
\end{proposition}
\bprf
Let $A_n$ and $\varphi_n=\Ad u_n\circ \xi^{(n)}:A_n \rightarrow A_{n+1} $ be as in the proof of Lemma \ref{lem:ifc} (i.e., $\xi^{(n)} = \xi_{d_n,m_n, r_{n,0},r_{n,1}, 1/2}\otimes 1_{c_{n+1}} $), $\widetilde{\varphi}_n$ the canonical map from $A_n$ into $A$, and $C_n= C([0,1])\otimes M_{d_n}.$ Remark that $\varphi_n$ was defined as a $*$-homomorphism $C_n \rightarrow C_{n+1}$ in the proof of Lemma \ref{lem:igh}. 
Let $(f_k)$ be a sequence of continuous functions in $C([0,1])^1_+$ such that $f_k(0) = f_k(1)=0 $ and  $f_k \nearrow \chi((0,1))$, where $\chi( (0,1))$ means the characteristic function on $(0,1)$.
We define a map $\widetilde{\psi}_n:C_n \rightarrow \pi_{\tau}(A)'' $ by 
\[ \widetilde{\psi}_n(g)= \lim_{k\rightarrow \infty }\pi_{\tau} \circ \widetilde{\varphi}_n ((f_k\otimes1_{d_n})\cdot g), \quad g \in C_n, \]
where the limit is taken in the strong operator topology. It follows that $\widetilde{\psi}_n $ is independent of the choice of $f_k$ and is an injective $*$-homomorphism, as follows. Since for $g \in {C_n}_+$, $\pi_{\tau}\circ \widetilde{\varphi}_n((f_k \otimes 1_{d_n})\cdot g) $ is a bounded increasing sequence of $\pi_{\tau}(A)$, the convergence follows for $g \in {C_n}_+$ and then for any $g \in {C_n}$. For $f_k$, $f_k' \in C([0,1])^1_+$ with the conditions $f_k(i)=f_k'(i)=0$, $i=0,1$, $f_k \nearrow \chi((0,1))$, and $f_k' \nearrow \chi((0,1))$, we have that
\begin{eqnarray}
\| \pi_{\tau} \circ \widetilde{\varphi}_n(((f_k-f_k')\otimes 1_{d_n})\cdot g) \|_2^2 &=&
\tau\circ \widetilde{\varphi}_n(((f_k-f_k')^2\otimes 1_{d_n})\cdot g^*g)\nonumber \\
&\leq& \| g \|^2 \tau\circ \widetilde{\varphi}_n ((f_k-f_k')^2
 \otimes 1_{d_n})\rightarrow 0 , \nonumber 
\end{eqnarray}
where $\|x\|_2^2=\tau(x^*x)$. Thus $\widetilde{\psi}_n$ is independent of $f_k$.
It is trivial that $\widetilde{\psi}_n$ is a linear map preserving $*$. For $f$, $g \in C_n$, we have that 
\begin{eqnarray}
\widetilde{\psi}_n(fg) &=& \lim \pi_{\tau} \circ \widetilde{\varphi}_n((f_k\otimes 1_{d_n})\cdot fg) \nonumber \\
&=&  \lim \pi_{\tau} \circ \widetilde{\varphi}_n((f_k^{1/2}\otimes 1_{d_n})\cdot f) \lim \pi_{\tau} \circ \widetilde{\varphi}_n((f_k^{1/2}\otimes 1_{d_n})\cdot g) \nonumber \\ &=& \widetilde{\psi}_n(f)\widetilde{\psi}_n(g). \nonumber
\end{eqnarray}
Thus $\widetilde{\psi}_n$ is a unital *-homomorphism. 
Because of the fact that the probability measure on $[0,1]$ obtained by $\tau\circ \widetilde{\varphi}_n$ is supported in $(0,1)$, we have that $\pi_{\tau}(\widetilde{\varphi}_n(f_k \otimes 1_{d_n})) \rightarrow 1 $ and that for any $f\in A_n$
\[\widetilde{\psi}_n (f) = \lim_{k\rightarrow \infty} \pi_{\tau} \circ \widetilde{\varphi}_n(f_k\otimes 1_{d_n})\cdot \pi_{\tau}\circ \widetilde{\varphi}_n(f) =\pi_{\tau} \circ \widetilde{\varphi}_n(f),\]
i.e., $\widetilde{\psi}_n|_{A_n}=\pi_{\tau} \circ \widetilde{\varphi}_n $.
Since 
\[ \|\widetilde{\psi}_n (g) \|_2^2 = \int \Tr_{d_n}(g(t)^* g(t)) d\mu (t),\]
for $g \in C_n$, where $\mu$ is the probability measure on $[0,1]$ obtained for $\tau |_{A_n} $, it follows that $\widetilde{\psi}_n(g)=0$ if and only if $g=0$ as $\mu $ is faithful, i.e., $\widetilde{\psi}_n$ is injective.

We assert that $\widetilde{\psi}_n = \widetilde{\psi}_{n+1} \circ \varphi_n$ on $C_n$. If a bounded sequence $g_k$ in $C_n$ converges to $g\in C_n$ pointwise on $(0,1)$, i.e., in the sense that $g_k(t) \rightarrow g (t)$ for $t\in (0,1)$, then it follows that $\widetilde{\psi}_n(g_k) $ converges to $\widetilde{\psi}_n(g) $ in the strong operator topology as 
\[ \|\widetilde{\psi}_n (g_k) - \widetilde{\psi}(g)  \|^2_2 \leq \int_{[0,1]} \|g_k(t) -g(t) \|^2 d \mu (t), \]
where $\mu$ is the probability measure as above. 
 Since $\widetilde{\psi}_n = \widetilde{\psi}_{n+1} \circ \varphi_n$ on 
$A_n$, it suffices to show that $A_n$ is dense in $C_n$ and $\varphi_n : C_n \rightarrow C_{n+1} $ is continuous in the pointwise topology on $(0,1)$. This follows easily. 
Thus we have that $\widetilde{\psi}_n=\widetilde{\psi}_{n+1}\circ \varphi_n$ and then we conclude that $\{ \widetilde{\psi}_n(C_n) \} $ is an increasing sequence of $C^*$-algebras satisfying 
\[ \pi_{\tau}(A) \subset \overline{(\bigcup \widetilde{\psi}_n(C_n))}^{\|\cdot\|}\subset \pi_{\tau} (A)''.\]

We define $B$ as $\overline{(\bigcup \widetilde{\psi}_n(C_n))}^{\|\cdot\|}$ and define a sequence of finite dimensional $C^*$-subalgebras $B_n$ of $B$ by 
\[B_n=\widetilde{\psi}_n\circ\Ad w_n(1_{C([0,1])}\otimes M_{d_n}), \]
where $w_n=u_{n-1}\xi^{(n-1)}(w_{n-1})\in C_n$, $n\in \N$, $w_0=u_0=1_{C_1}$, and $\xi^{(0)} = \id _{C_1} $.
Since \[ \varphi_{N,n}(f)= \Ad w_{N+1} \circ \xi^{(N)}\circ \xi^{(N-1)}\circ \cdots \circ\xi^{(n)}\circ \Ad w_n^*(f),\]
for any $f \in C_n$, where $\varphi_{N,n}=\varphi_N\circ \varphi_{N-1}\circ \cdots \circ \varphi_n ,$ and since  
\[ B_n=\widetilde{\psi}_{n+1}\circ\Ad u_n \circ \xi^{(n)} \circ \Ad w_n (1\otimes M_{d_n})= \widetilde{\psi}_{n+1}\circ \Ad w_{n+1} \circ \xi^{(n)}(1\otimes M_{d_n}),\] we have that $B_n\subset B_{n+1}$ and $1_{B_n}=1_{B_{n+1}}$.
By the definition of $\xi^{(n)}$, for $f\in C_n$ there exists a large natural number $m$ such that $\xi^{(m)}\circ\cdots \circ \xi^{(n)}(f)$ is almost contained in $1_{C([0,1])}\otimes M_{d_{m+1}}$  (i.e., there exists $ x \in 1_{C[0,1]}\otimes M_{d_{m+1}}$ such that $\|\xi^{(m)}\circ\cdots \circ \xi^{(n)}(f)-x\|<\varepsilon  $). Then 
\[ \widetilde{\psi}_n(f)=\widetilde{\psi}_{m+1}\circ\varphi_{m,n}(f)=\widetilde{\psi}_{m+1}\circ \Ad w_{m+1} \circ \xi^{(m)}\circ \xi^{(m-1)}\circ \cdots \circ\xi^{(n)}\circ \Ad w_n^* (f)\] is almost contained in $B_{m+1}$. It follows that $\bigcup B_n$ is a norm-dense subalgebra of $B$. Since $B_n \cong M_{d_n}$, $B$ is a UHF algebra.

 We now ensure the last condition $(\pi_{\tau}(A) \cap B_n')''= \pi_{\tau}(A)''\cap B_n'$, where the part ``$\subset$" is obvious. For $x\in \widetilde{\psi}_m(C_m)\cap B_n' $,  
there exists $y\in C_m\cap \varphi_{m-1,n}\circ \Ad w_n(1\otimes M_{d_n})' $ such that $x=\widetilde{\psi}_m(y)$. Since $\pi_{\tau}\circ \widetilde{\varphi }_m( (f_k\otimes1_{d_m})\cdot y) $ is contained in $\pi_{\tau}(A) \cap B_n'$, it follows that $x \in (\pi_{\tau}(A)\cap B_n')''$. Since for any $x \in \pi_{\tau}(A)''\cap B_n'$ there exists $x_m \in \widetilde{\psi}_m(C_m)\cap B_n'$ such that $x_m \rightarrow x\ ({strongly})$, we have that $\pi_{\tau}(A)''\cap B_n'\subset (\pi_{\tau}(A)\cap B_n' )''$.\eprf

The following stability lemma was proved by Connes in \cite{Co}.

\begin{lemma}\label{lem:Cs}
Let $M$ be the injective type $I\hspace{-.1em}I_1$-factor and let $\tau$ be the unique tracial state of $M$. Suppose that an automorphism $\theta$ of $M$ is aperiodic (i.e., $\theta^n$ is outer automorphism for all $n \in \N$ ). Then for $\varepsilon >0$ and a finite subset $F$ of $M$ there exist $\delta>0$ and a finite subset $G$ of $M$ satisfying the following condition: for any unitary $u\in M $ with $\|[u,y]\|_2 < \delta,\ y\in G$ there exists a unitary $v
\in M$ such that $\|u- v\theta(v^*) \|_2 < \varepsilon$ and $\| [v, x] \|_2 <
\varepsilon ,\ x\in F $, where $\|x\|_2 = \tau(x^*x)^{1/2} $. Furthermore, if $F$ is empty then it may be assumed that $G$ is empty.   
 
\end{lemma}

In the proof of the main theorem, we need the following fundamental lemma for the hyperfinite type $I\hspace{-.1em}I_1$-factor.

\begin{lemma}\label{lem:hyperfinite}
Let $M$ be the hyperfinite type $I\hspace{-.1em}I_1$-factor, $\tau$ the unique tracial state of $M$, and $\theta$ an automorphism of $M$. 
For a finite subset $F$ of $M$ and $\varepsilon >0$, there exists a unitary $u$ in $M$ such that $\|\theta (f) - \Ad u(f)\|_2 < \varepsilon $ for any $f\in F$ .
\end{lemma}
\bprf
Let $F$ be a finite subset of $M$ and $\varepsilon>0$. 
There exists a finite dimensional $C^*$-subalgebra $B$ of $M$ such that
$B$ is isomorphic to a full matrix algebra and for any $f\in F$ there exists $g_f\in B$ such that $\|f - g_f\|_2 <2^{-1}\varepsilon $.
Let $\{e_{i,j}\}$ be a system of matrix units of $B$. Since $\tau(e_{1,1})=\tau
\circ \theta (e_{1,1})$, there exists $v\in M$ such that $v^*v=e_{1,1}$, $vv^*=
\theta (e_{1,1})$. We define a unitary $u\in M$ by $\displaystyle u=\sum_i \theta (e_{i,1})v e_{1,i}$. Then we have that $\theta (b) = \Ad u (b)$ for any $b\in B$. Thus we have 
\[\|\theta (f) - \Ad u (f) \|_2 \leq \| \theta(f-g_f)\|_2+ \| \Ad u(f-g_f) \|_2 <
\varepsilon .\] \eprf

The proof of the following main theorem is based on the intertwining argument for automorphisms of AF-algebras which is introduced by Evans and Kishimoto (\cite{EK}). We recall that $\WInn (A)$ denotes the normal subgroup of $\Aut (A)$ consisting of $\alpha \in \Aut (A) $ with $\pi_{\tau}\alpha = \Ad u \circ \pi_{\tau} $ for some $u \in U(\pi_{\tau} (A)'') $.

\begin{theorem}\label{thm:wa}
Let $A$ be a unital simple projectionless $C^*$-algebra which is the inductive limit of a sequence of dimension drop algebras, as in Lemma \ref{lem:ifc}, such that $K_1(A)$ is the inductive limit of a sequence of finite cyclic groups and $A$ has a unique tracial state $\tau$. 
If $\alpha$, $\beta \in \Aut(A)$ are aperiodic in the quotient $\Aut(A)/\WInn(A)$, then for any $\varepsilon >0$ there exist an approximately inner automorphism $\sigma \in \Aut (A)$ and $\gamma \in \WInn (A) $ with $W \in U(\pi_{\tau}(A)''
)$ such that $\pi_{\tau} \gamma =\Ad W \circ \pi_{\tau}$, $\|W-1 \|_2< \varepsilon $, and
\[\alpha= \gamma \circ \sigma \circ \beta \circ \sigma^{-1}. \]
\end{theorem}
\bprf
By Proposition \ref{prop:eUHF} there exists a UHF $C^*$-subalgebra $B$ of $\pi_{\tau}(A)''$ and an increasing sequence of matrix algebras $B_n \subset B$ such that $\pi_{\tau} (A) \subset B $, $1_{B_n}=1_{B_{n+1}}$, $\displaystyle B= \overline{(\bigcup B_n)}^{\|\cdot\|} $, and $(\pi_{\tau}(A) \cap B_n)''= \pi_{\tau}(A)'' \cap B_n'$. Let $\{e_{i,j}^{(n)}\}_{i,j}$ be a system of matrix units of $B_n$ and  $\varepsilon >0 $, and set $B_0=\C 1$. 

In the following we will identify $\pi_{\tau} (A) $ with $A$. Hence if $\gamma \in \Aut (A) $, $\overline{\gamma} $ denotes the weak extension of $\gamma $ to an automorphism of $A''$ ($=\pi_{\tau} (A)'' $). We denote by $x \in_{\delta} B_k $ the condition that there is $y\in B_k$ with $\|x-y\| < \delta $.

We shall construct, inductively, finite subsets $ F_{2l+1}, G_{2l+2} $ of $ A''$,  $u_{2l+1}$, $u_{2l+2}$, $v_{2l+1}$, $v_{2l+2} \in U(A) $,  $k_{2l+1},k_{2l+2} \in \Z_+$, and $\delta_{2l+1}, \delta_{2l+2}>0$ for $l=0,1,2,...$ satisfying the following conditions: 
Setting $G_0=\{1\}$, $v_0=1$, $k_0=0$, $\delta_0=\varepsilon$, 
$\alpha_{-1} = \alpha$, $\beta_0= \beta $, $w'_0=w'_1=w'_2=1$, for $l \geq0$
\[ \alpha_{2l+1}= \Ad u_{2l+1}\circ \alpha _{2l-1},\quad \beta_{2l+2}= \Ad u_{2l+2} \circ \beta_{2l}, \] 
\[\sigma^{(0)}_l= \Ad v_{2l}v_{2l-2}\cdots v_0,\quad \sigma^{(1)}_l=\Ad v_{2l+1}v_{2l-1}\cdots v_1,\] 
\[w_{2l+1}=u_{2l+1}\alpha_{2l-1}(v_{2l+1})v_{2l+1}^* \in A, \quad w_{2l+2}=u_{2l+2}\beta_{2l}(v_{2l+2})v_{2l+2}^* \in A,\]
and for $l \geq 1 $, $i=0,1$
\[w_{2l+1+i}'=w_{2l+1+i}\Ad v_{2l+1+i} (w_{2l-1+i}')\in A,\] 
the conditions are given by 
\begin{description}
\item[(1)]\label{(1)} $G_{2l} \subset F_{2l+1},\quad \{e_{i,j}^{(k_{2l+1})}\} \subset F_{2l+1}$,  
\item[(1)']\label{(1)'} $F_{2l+1}\subset G_{2l+2},\quad \{e_{i,j}^{(k_{2l+2})}\} \subset G_{2l+2},$
\item[(2)]\label{(2)} $\| \Ad u_{2l+1}\circ \overline{\alpha}_{2l-1}(x) - \overline{\beta}_{2l}(x) \|_2<2^{-1}\delta_{2l+1}$ for any $ x\in F_{2l+1}$,
\item[(2)']\label{(2)'} $\| \overline{\alpha}_{2l+1} (x) - \Ad u_{2l+2} \circ \overline{\beta}_{2l}(x) \|_2 < 2^{-1}\delta_{2l+2} $ for any $x\in G_{2l+2}$, 
\item[(3)]\label{(3)} $\|u_{2l+1}-v_{2l+1}\alpha_{2l-1}(v_{2l+1}^*)\|_2 <2^{-(l+2)}\varepsilon,\quad v_{2l+1} \in A \cap B_{k_{2l}}'$,
\item[(3)']\label{(3)'} $\|u_{2l+2}-v_{2l+2}\beta_{2l}(v_{2l+2}^*)\|_2 <2^{-(l+3)}\varepsilon,\quad v_{2l+2} \in A \cap B_{k_{2l+1}}'$,
\item[(4)]\label{(4)} $\sigma^{(0)}_l(e_{i,j}^{(m)}) \in_{2^{-(l+2)}} B_{k_{2l+1}}$ for any $0\leq m\leq k_{2l}$, 
$\\ w_{2l}'\in_{2^{-(l+3)}\varepsilon} B_{k_{2l+1}}\quad, \quad k_{2l+1}\geq 2k_{2l},$
\item[(4)']\label{(4)'} $\sigma^{(1)}_l(e_{i,j}^{(m)}) \in_{2^{-(l+3)}} B_{k_{2l+2}} $ for any $0\leq m\leq k_{2l+1}$, 
$\\ w_{2l+1}'\in_{2^{-(l+4)}\varepsilon} B_{k_{2l+2}} $, $\quad k_{2l+2}\geq 2k_{2l+1}$,
\item[(5)]\label{(5)} $\delta_{2l+1}\leq 2^{-1} \delta_{2l}$, and
 
 if $u\in U(A'')$ satisfies that $\|[u,x]\|_2< \delta_{2l+1}$ for any $ x \in \overline{\beta}_{2l}(F_{2l+1})$,   
then there exists $v'\in U(A'') $ such that 
 \[\|[v',e_{i,j}^{(k_{2l+1})}]\|_2 < 2^{-(l+7)} (\dim(B_{k_{2l+1}}))^{-1/2}\varepsilon,\quad \| u-v'\overline{\beta}_{2l}(v'^*) \|_2 < 2^{-(l+4)}\varepsilon,\] 

\item[(5)']\label{(5)'} $\delta_{2l+2}\leq 2^{-1}\delta_{2l+1}$, and
 
if $u\in U(A'')$ satisfies that $\|[u,x]\|_2< \delta_{2l+2}$ for any $ x \in \overline{\alpha}_{2l+1}(G_{2l+2})$, 
then there exists $v'\in U(A'') $ such that
\[\|[v',e_{i,j}^{(k_{2l+2})}]\|_2 < 2^{-(l+7)} (\dim(B_{k_{2l+2}}))^{-1/2}\varepsilon,\quad \| u-v'\overline{\alpha}_{2l+1}(v'^*) \|_2 < 2^{-(l+4)}\varepsilon.\] 
\end{description} 
For $l=0$, set 
$ F_1=\{1 \},\ G_2=\{1 \},\quad u_1=u_2=1_{A},\quad v_1=v_2=1_{A},\quad k_1=k_2=0,\quad \delta_1=2^{-1}\varepsilon$, and $\delta_2=2^{-2}\varepsilon$. The conditions (1)$\sim$(4)' are trivially satisfied for $l=0$ and (5) and (5)' follows from Lemma \ref{lem:Cs} since $\overline{\alpha}$, $\overline{\beta}$ are aperiodic automorphisms of $ A''$. Assuming that we have constructed 
\[F_{2l+1},G_{2l+2},u_{2l+1},u_{2l+2},v_{2l+1},v_{2l+2},k_{2l+1},k_{2l+2},\delta_{2l+1},\delta_{2l+2} \]
 for $l\leq n-1$ which satisfy (1)$\sim$(5)', we proceed as follows:
Note that $u_m$, $v_m$, $k_m$, $\delta_m$, $w_m$, and $w_m'$ are now given for $m\leq 2n$. Since $\sigma _n ^{(0)} $ leaves $B$ invariant and $w_{2n}'\in B$, we can find a $k_{2n+1}$ satisfying (4) for $l=n$, i.e., $ k_{2n+1} \geq 2k_{2n} $ and for any $ m\leq k_{2n}$, 
\[\sigma_n^{(0)} ( e_{i,j}^{(m)} ) \in_{2^{-(n+2)}} B_{k_{2n+1}},\quad w_{2n}' \in_{2^{-(n+3)}\varepsilon} B_{k_{2n+1}}.\]
 Since $\overline{\beta}_{2n} $ is an aperiodic automorphism of $A'' $, there exists a finite subset $F_{2n+1}$ of $A''$ and $\delta_{2n+1}>0 $ which satisfy the conditions (1), (5) for $l=n$, by Lemma \ref{lem:Cs}.
Since $A''$ is the hyperfinite type $I\hspace{-.1em}I_1$-factor, by Lemma \ref{lem:hyperfinite} there is a unitary $u_{2n+1}' \in A'' $ such that for any $ x\in F_{2n+1} $, 
\[ \|\Ad u_{2n+1}'\circ \overline{\alpha}_{2n-1}(x) - \overline{\beta}_{2n}(x) \|_2 <2^{-1}\delta_{2n+1}. \] 
Then there exists a unitary $u_{2n+1} \in A$ satisfying (2) for $l=n$ by the Kaplansky's density theorem. By the assumption (2)' for $l=n-1$ and (1), (2) for $l=n$, $u_{2n+1} $ satisfies that for any $x \in  G_{2n}$,
\begin{eqnarray}
\|[u_{2n+1},\overline{\alpha}_{2n-1}(x)]\|_2 &=& \|\Ad u_{2n+1}\circ \overline{\alpha}_{2n-1}(x)-\overline{\alpha}_{2n-1}(x)\|_2  \nonumber \\
&<& \|\Ad u_{2n+1}\circ\overline{\alpha}_{2n-1}(x)-\overline{\beta}_{2n}(x)\|_2 + 2^{-1}\delta_{2n} \nonumber \\
&<& 2^{-1}\delta_{2n+1} + 2^{-1}\delta_{2n}<\delta_{2n}. \nonumber
\end{eqnarray} 
Then by (5)' for $l=n-1 $ there exists a unitary $v'\in U( A'')$ satisfying the condition in (5)' for $l=n-1$ with $u=u_{2n+1}$. 
Let $\Phi_{B_{k_{2n}}}: A'' \rightarrow  A'' \cap B_{k_{2n}}'$ be the conditional expectation. Since $v'$  almost commutes with $\{
e_{i,j}^{(k_{2n})} \} $ in the sense that 
\[ \|[v',e_{i,j}^{(k_{2n})}]\|_2 < 2^{-(n+6)} (\dim(B_{k_{2n}}))^{-1/2}\varepsilon,\] 
we have that  
\[\|\Phi_{B_{k_{2n}}} (v') - v' \|_2\leq \dim(B_{k_{2n}})^{-{1/2}}\sum_{i,j}\|[v',
e_{i,j}^{(k_{2n})}]\|_2< 2^{-(n+6)}\varepsilon. \]
By the polar decomposition of $\Phi_{B_{k_{2n}}} (v')$ we obtain a unitary $v'' \in U(A '' \cap B_{k_{2n}}' )$ such that $ v''\cdot |\Phi_{B_{k_{2n}}}(v')|=\Phi_{B_{k_{2n}}} (v')$. Since 
\[ \|\Phi_{B_{k_{2n}}}(v')^*\Phi_{B_{k_{2n}}}(v')-1\|_2<2^{-(n+5)}\varepsilon,\]
we have that $\||\Phi_{B_{k_{2n}}}(v')|-1\|_2<2^{-(n+5)}\varepsilon$.  
Thus $v''$ satisfies that
\[\|v''-v'\|_2<\|\Phi_{B_{k_{2n}}}(v')-v'\|_2+2^{-(n+5)}\varepsilon<2^{-(n+4)}\varepsilon.\]
By $ A'' \cap B_{k_{2n}}' = (A \cap B_{k_{2n}}')''$,
 there exists a unitary $v_{2n+1} \in U( A \cap B_{k_{2n}}') $ such that $ \| v_{2n+1} - v' \|_2 < 2^{-(n+4)}\varepsilon $. Then we have that 
\[ \|v_{2n+1}\alpha_{2n-1}(v_{2n+1}^*)-u_{2n+1}\|_2<\|v'\overline{\alpha}_{2n-1}(v'^*)-u_{2n+1}\|_2+2^{-(n+3)}\varepsilon<2^{-(n+2)}\varepsilon,\]
 which is the condition (3) for $l=n$.
We have now constructed $F_m$, $u_m$, $v_m$, $k_m$, and $\delta_m$ for $m= 2n+1$ satisfying the conditions with no primes for $l=n$ and will construct $G_m$, $u_m$, $v_m$, $k_m$, and $\delta_m$ for $m=2n+2$ satisfying the conditions (1)' $\sim$ (5)' for $l=n$ in a similar way. Since $\sigma^{(1)}_n $ leaves $B$ invariant and $w_{2n+1}'\in B$,   we can find a large natural number $k_{2n+2}\ \geq 2k_{2n+1} $ such that for any $ m\leq k_{2n+1}$,
\[\sigma_n^{(1)} ( e_{i,j}^{(m)} ) \in_{2^{-(n+3)}} B_{k_{2n+2}},\quad w_{2n+1}' \in_{2^{-(n+4)}\varepsilon} B_{k_{2n+2}}.\]
By Lemma \ref{lem:Cs} there are a finite subset $G_{2n+2}$ of $A''$
and $\delta_{2n+2}>0 $ which satisfy (1)', (5)' for $l=n$. By the same reasoning as above there exists a unitary $u_{2n+2} \in A$ satisfying (2)' for $l=n$. By the condition (2) for $l=n$ and (1)', (2)' for $l=n$, it follows that  
\[ \|[u_{2n+2},x]\|_2 < 2^{-1}\delta_{2n+2} + 2^{-1}\delta_{2n+1} <
\delta_{2n+1},\quad x \in \overline{\beta}_{2n}( F_{2n+1})  .\]
By (5) for $l=n $ and by $ A'' \cap B_{k_{2n+1}}' = (A \cap B_{k_{2n+1}}')''$, there exists a unitary $v_{2n+2} \in U( A \cap B_{k_{2n+1}}') $ with the condition of (3)' for $l=n$. This completes the induction.

We note that the linear span of the set $\bigcup_{n} F_{2n+1} $ is strongly dense in $A''$. By the conditions (4), (4)' for $l=n$, we have that for $i=0,1$, $m\leq 2n+i$, and $1\leq p,q \leq \dim (B_{k_{m}})^{1/2}$, 
there exists $f_{p,q}^{(m)}\in B_{k_{2n+1+i}}$ such that  
\[\|f_{p,q}^{(m)}-\sigma_n^{(i)}(e_{p,q}^{(m)})\|<2^{-(n+2+i)}.\]
By (3), (3)', we have that
\begin{eqnarray} 
\|\sigma^{(i)}_{n+1}(e_{p,q}^{(m)})-\sigma^{(i)}_n(e_{p,q}^{(m)})\|&=&\|\Ad v_{2n+2+i}\circ \sigma^{(i)}_{n}(e_{p,q}^{(m)})-\sigma^{(i)}_n(e_{p,q}^{(m)})\| \nonumber \\
&\leq& \|\Ad v_{2n+2+i} (f_{p,q}^{(m)}) - f_{p,q}^{(m)}\|+2^{-(n+1+i)}\nonumber\\
&=&2^{-(n+1+i)}.\nonumber
\end{eqnarray} 
 Thus for $i=0,1$ and $x\in\bigcup_m B_{m}$, $\sigma_n^{(i)}(x)$ is a norm Cauchy sequence in $B$. Then we can define endomorphisms $\sigma_i,\ i=0,1,$ of $B$ by 
\[\sigma_i(x)= \lim_{n} \sigma^{(i)}_n (x)\quad({\rm norm}),\ x\in B. \] 
It follows that $\sigma_i(A) \subset A,\ i=0,1,$ since $v_n \in A $. By (3), (3)', we can also define endomorphisms $\sigma_i^{-1}$  of $B$ by 
\[\sigma_i^{-1} (x)= \lim_{n} (\sigma^{(i)}_n)^{-1}(x),\quad x\in B.\] We can easily prove that $\sigma_i\circ \sigma_i^{-1}=id_B= \sigma_i^{-1} \circ \sigma_i$ and thus $\sigma_i^{-1}$ is indeed the inverse of $\sigma_i$. Since we also have that $\sigma_i^{-1}(A) \subset A$, we conclude that $\sigma_i$ is an automorphism of $B$ and restricts to an automorphism of $A$.

By (3), (3)', we have that $\| w_{2n+1+i}-1\|_2 < 2^{-(n+2+i)}\varepsilon $ for $i=0,1$ and by (4)' for $l=n-1$ and by (4) for $l=n$, there exists $w_{2n-1+i}''\in B_{k_{2n+i}}$ for $i=0,1$ such that $\|w_{2n-1+i}''-w_{2n-1+i}'\|<2^{-n+3}\varepsilon.$
We have that for $n\in \N$ and $i=0,1$, 
\begin{eqnarray}
\|w_{2n-1+i}'-w_{2n+1+i}'\|_2 &=& \|w_{2n-1+i}'-w_{2n+1+i}\Ad v_{2n+1+i}(w_{2n-1+i}')\|_2  \nonumber \\
&<& \|w_{2n-1+i}'-\Ad v_{2n+1+i}(w_{2n-1+i}')\|_2+2^{-(n+2+i)}\varepsilon \nonumber\\ 
&<& 2^{-(n+2)}\varepsilon + 2^{-(n+2+i)}\varepsilon \leq 2^{-(n+1)}\varepsilon  \nonumber
\end{eqnarray} 
Thus we obtain unitaries $\widetilde{w}_i \in U(A'') $ for $i=0,1$ by 
\[\widetilde{w}_i= \lim_{n} w_{2n+1+i}'\quad ({\rm strongly}),\] 
which satisfies
\[\|1-\widetilde{w}_i\|_2 \leq \sum_{n=1}^{\infty}\|w_{2n-1+i}' - w_{2n+1+i}'\|_2
 <2^{-1}\varepsilon, \quad i=0,1. \]

By the definition of $w_{2n+1+i}'$, we have that \[\overline{\alpha}_{2n+1}=\Ad w_{2n+1}' \circ \overline{\sigma}_{n}^{(1)} \circ \overline{\alpha} \circ (\overline{\sigma}_{n}^{(1)})^{-1}, \quad \overline{\beta}_{2n} = \Ad w_{2n}'\circ \overline{\sigma}_{n}^{(0)} \circ 
\overline{\beta} \circ (\overline{\sigma}_{n}^{(0)})^{-1},\]
and by (2), (2)', we have that for any $n\in \N$ and $x\in F_{2n+1}$,
\[\|\Ad w_{2n+1}' \circ \overline{\sigma}_{n}^{(1)} \circ \overline{\alpha} \circ (\overline{\sigma}_{n}^{(1)})^{-1}(x)- \Ad w_{2n}'\circ \overline{\sigma}_{n}^{(0)} \circ \overline{\beta} \circ (\overline{\sigma}_{n}^{(0)})^{-1}(x)\|_2 < 2^{-1}\delta_{2n+1}.\]
Since $\delta_{2n+1} \rightarrow 0$ by $\delta_{2n+1}\leq 2^{-1} \delta_{2n-1}$, we have that for any $x\in A''$, 
\[ \Ad \widetilde{w}_1 \circ \overline{\sigma}_1 \circ \overline{\alpha} \circ (\overline{\sigma}_1)^{-1}(x) = \Ad \widetilde{w}_0 \circ \overline{\sigma}_0 \circ \overline{\beta} \circ (\overline{\sigma}_0)^{-1}(x),\]
 which implies that on $A''$ 
\[ \overline{\alpha} = \Ad (\overline{\sigma}_1)^{-1}(\widetilde{w}_1^*\widetilde{w}_0) \circ (\overline{\sigma}_1)^{-1}\circ \overline{\sigma}_0\circ \overline{\beta} \circ (\overline{\sigma_0})^{-1}\circ \overline{\sigma_1}. \] 
Set $\sigma=\sigma_{1}^{-1}|_{A}\circ\sigma_{0}|_{A} \in \Aut (A)$ and $W=(\overline{\sigma}_1)^{-1}(\widetilde{w}_1^*\widetilde{w}_0) \in U(A'')$. Then we have that $\Ad W(A)=A$, $\|W-1\|_2<\varepsilon$, and $\alpha=\Ad W\circ\sigma\circ\beta\circ \sigma^{-1}$.  This completes the proof. \eprf

\section{The Rohlin property and a certain class of TAF-algebras}\label{sec:P}
First, we recall the definition of the Rohlin property.
 
\begin{definition}\label{def:R}
Let $A$ be a unital $C^*$-algebra and $\alpha$ be an automorphism of $A$. We say 
that $\alpha$ has the Rohlin property if there exists $k \in \N $  which satisfies the following condition: For any finite subset $F$ of $A$ and any $\varepsilon >$ 0 there exists projections $\{ e_j^{(0)}\ :j = 0,...,k-1 \} \cup \{ e_j^{(1)}\ :j=\ 0,...,k \}$ in $A$ such that
\[\sum_{j=0}^{k-1}e_j^{(0)}+\sum_{j=0}^k e_j^{(1)}=1_A ,\]
\[\| \alpha(e_j^{(i)})-e_{j+1}^{(i)} \| < \varepsilon, \quad i=0,1, \quad j=0,...,k+i-2,\]
\[\| [a\ , e_j^{(i)}]\| <\varepsilon\]
for all $a\in F$, $i=0,1$ and $j=0,1,...,k+i-1$.

These $\{e_j^{(0)}\}$ and $\{e_j^{(1)}\}$ are called Rohlin towers. This definition of 
the Rohlin property was defined by A.Kishimoto in \cite{Kis2} and he showed the 
stability for automorphisms of a unital simple $A\T$-algebra of real rank zero 
with the Rohlin property defined as above.
\end{definition}

\begin{example}\label{exR} Let $A$ be the UHF-algebra $\bigotimes_{p:\text{prime}} M _p $ which corresponds to the supernatural number obtained by products of all prime numbers. Define an automorphism $\alpha $ of $A$ by \[ \alpha:=\bigotimes_{p:\text{prime}}\ \Ad s_p , \]
where $s_p$ is the shift operator in $M_p$.
Then it follows that $\alpha $ has the Rohlin property. But this example does not satisfy a version of Rohlin property defined through one tower. More information about the Rohlin property appears in \cite{Izu}.
\end{example}

In order to prove the main result in the next section, we need the following property (SI), where SI means small isometry. For unital simple $A\T$-algebras of real rank zero with a unique tracial state, we can see a similar property in Lemma 4.6 \cite{Kis2}.  

\begin{definition}\label{def:P}
Let $A$ be a separable unital $C^*$-algebra. We say that the {\it order on projections over $A$ is determined by traces } if $p, q \in P(M_{\infty}(A))$ with $\tau(p)<\tau(q)$ for any $\tau \in T(A) $ satisfy $p\preceq q$ (i.e., there exists $w \in A$ such that $w^* w= p$ and $w w^*\leq q $), see \cite{Bl} and \cite{OP1}. We say that $A$ has the {\it property (SI)}, if $A$ satisfies that: if a unital $C^*$-algebra $B$ of real rank zero is given such that there exists a unital embedding from $A$ into $B$ and the order on projections over $B$ is determined by traces, and any sequences $(e_n)_n$, $(f_n)_n \in P( B^{\infty} \cap A')$ with $e_n$, $f_n \in B_{\rm sa}$  satisfy the following conditions:
\begin{description}
\item[(1)]\label{SI:1} there exists $c>0$  satisfying that: for any $p \in P(A) $ and $\varepsilon>0$ there exists $N_p\in\N$ such that $\ c\tau(p) - \varepsilon \leq \tau(pf_np)$ for any $n \geq N_p $ and $\tau \in T(B)$,
\item[(2)]\label{SI:2} $\tau (e_n)  \rightarrow 0$ for any $\tau \in T(B) $,  
\end{description} 
then there exists a $(w_n)_n \in  B^{\infty} \cap A'$
such that 

\[ (w_n)^*(w_n)=(e_{n}),  \quad (w_n)(w_n)^*\leq (f_{n}),\quad (w_n)^2=0.\]
\end{definition}

\begin{remark}\label{rem:SI} The above definition is equivalent to the following seemingly weaker condition: For any unital $C^*$-algebra $B$ of real rank zero and any $(e_n)$, $(f_n)\in P( B^{\infty} \cap A')$ with the same conditions in the above definition, there exist a sequence $(m_n)$ in $\N$ and $(w_n)_n\in B^{\infty} \cap A'$ such that  
\[ (w_n)^*(w_n)=(e_{m_n}),  \quad (w_n)(w_n)^*\leq (f_{m_n}),\quad (w_n)^2=0.\]

This equivalence follows from the following reason. Assume that $A$ does not have the property (SI). Then there exists a unital $C^*$-algebra $B$ of real rank zero, $(e_n)$, $(f_n)\in P(B^{\infty}\cap A')$ with the above conditions, finite subset $F\subset A^1$, and $\varepsilon_0> 0$ satisfying that: for any $n\in \N$ there exists $N \geq n$ such that any $w \in B$ satisfies $\|w^*w -e_N \| >\varepsilon_0$, $\|ww^*-ww^*f_N\| >\varepsilon_0 $, $\|w^2\|> \varepsilon_0$, or $\|[w,a]\|>\varepsilon_0$, for some $a \in F$. Thus there exists a sequence $(l_n)_n$ in $\N$ satisfying that: any subsequence $(m'_n)$ of $(l_n)$ and any $(w_n) \in B^{\infty} \cap A'$ satisfy that $(w_n)^*(w_n)\neq (e_{m'_n})$, $(w_n)(w_n)^* \nleq (f_{m'_n})$, or $(w_n)^2\neq 0$. But $(e_{l_n})_n$ and $(f_{l_n})_n$ also satisfy the conditions (1) and (2). Hence $A$ does not satisfy the latter condition.   
\end{remark}
The basic idea for the proof of the following lemma appears in the proof of Lemma 2.8 in \cite{OP2} or the proof of Theorem 4.5 in \cite{Kis}.    
\begin{lemma}\label{lem:Orthogonal}
Let $A$ be a unital $C^*$-algebra of real rank zero. 
Suppose that the order on projections over $A$ is determined by traces.
If non-zero projections $e$, $f$ in $A$ satisfy $2\tau(e) < \tau (f) $ for any
$\tau \in T(A)$, then for any $\varepsilon >0$ there exists a partial isometry $w\in A$ such that 
\[ w^*w =e,\quad ww^*\leq f, \quad \|w^2\| < \varepsilon. \]
\end{lemma}
\bprf
Let $e$ and $f$ be projections in $A$ as in the statement. 
We set $\varepsilon' = 18^{-1}\varepsilon^2$ and define continuous functions  $g$, $h\in C([0,1])_+^1$ by
\begin{eqnarray}
g(t)=\left\{ \begin{array}{ll}
t/{\varepsilon'},\quad & 0\leq t < \varepsilon' \\
1,\quad  & \varepsilon' \leq t< 1, \\
\end{array} \right.\nonumber
\end{eqnarray}   
\begin{eqnarray}
h(t)=\left\{ \begin{array}{ll}
0,\quad & 0\leq t < \varepsilon' \\
t/{\varepsilon'}-1,\quad  & \varepsilon'\leq t < 2\varepsilon' \\
1, &  2\varepsilon' \leq t \leq 1. \\
\end{array} \right.\nonumber
\end{eqnarray}
Because $\overline{h(fef)Ah(fef)}$ is a $C^*$-algebra of real rank zero, there exists $q\in P(\overline{h(fef)Ah(fef)})$ such that \[\|qh(fef) -h(fef)\| < \varepsilon/3.\]
Set $p=f-q\in P(fAf)$.
Since $\|h(fef)fe -fe \|^2$ 
\begin{eqnarray}
&=&\| h(fef)fefh(fef)-h(fef)fef-fefh(fef)+fef\| \nonumber \\
&\leq & 2 \|h(fef)fef- fef\| \leq 2 \varepsilon', \nonumber
\end{eqnarray}
 we have that $\|pe \| = \|fe- qfe\|$
\[ \leq \| qfe-qh(fef)fe\|+\|qh(fef)fe-h(fef)fe\| + \|h(fef)fe-fe\| <\varepsilon.\]

Since the order on projections over $A$ is determined by traces and 
\begin{eqnarray}
\tau(p) &=& \tau(f-q) = \tau(f)- \tau(g(fef)^{1/2}qg(fef)^{1/2}) \nonumber \\
&\geq& \tau(f) - \tau (g(fef)) = \tau(f) - \tau(g(efe)) \nonumber \\
&\geq& \tau(f) - \tau(e) > \tau(e) , \nonumber
\end{eqnarray}
there exists $w\in A$ such that $w^*w= e$ and $ww^* \leq p $. 
This $w$ satisfies that $\|w^2\| \leq \|ep\| < \varepsilon $.
\eprf
\begin{proposition}\label{Prop:AT}
Any unital simple A$\T$-algebra of real rank zero satisfies the property (SI).
\end{proposition} 
\bprf
Let $A$ be a unital simple infinite-dimensional A$\T$-algebra of real rank zero and $B$ a unital $C^*$-algebra of real rank zero such that the order on projections over $B$ is determined by traces and $A \subset_{\it unital} B$. Let $(e_n)_n$, $(f_n)_n  \in P(B^{\infty} \cap A')$, and $c>0$ satisfying the following conditions: for any $p\in  P(A)$ and any $\varepsilon > 0 $ there exists a large natural number $N_p$ such that $c\tau(p) - \varepsilon \leq \tau (p f_n p)$ for any $n \geq N_p$ and any $\tau \in T(B)$, and  $\tau(e_n) \rightarrow 0 $ for any $\tau \in T(B) $, where we have assumed that $e_n$, $f_n$ are all projections.
 Let $F$ be a finite subset of $A^1$ and $\varepsilon>0$. By Remark \ref{rem:SI} it suffices to show that, for any $n\in \N$, there exist $w \in B$ and $m\geq n$ satisfying the following conditions:  $\|[w,x]\| < \varepsilon $ for any $x \in F$, 
\[ \|w^*w- e_{m}\|<\varepsilon,\quad \|ww^*-ww^*f_m\|< \varepsilon,\quad \|w^2\| < \varepsilon. \]
Indeed, for an increasing sequence of finite subsets  $(F_n)_{n\in \N}$ of $A^1 $ such that $\bigcup F_n $ is a dense subset of $A^1 $ and for a decreasing sequence $(\varepsilon_n)_{n\in\N} $ such that $\varepsilon_n >0$ and $\varepsilon_n \rightarrow 0$, we can inductively obtain $w_n \in B$, and $m_n \in \N$ such that $m_n>m_{n-1}$, $\|{w_n}^*w_n - e_{m_n} \| < \varepsilon_n$, $\|w_n{w_n}^*- w_n{w_n}^*f_{m_n}\| <\varepsilon_n$, $\|w_n^2\|<\varepsilon_n $, and $\|[w_n, f]\| < \varepsilon_n $ for any $f\in F_n$. Thus we have that $(w_n)_n \in B^{\infty} \cap A'$, $(w_n)^*(w_n)=(e_{m_n})$, $(w_n)(w_n)^* \leq (f_{m_n})$, and $(w_n)^2=0$, which is equivalent to (SI) by Remark \ref{rem:SI}.  

Since $A$ is a unital simple infinite-dimensional A$\T$-algebra of real rank zero, there is  an increasing sequence $(A_n)_n$ of $C^*$-subalgebras $A_n \cong B_n \otimes C(\T) $ with finite-dimensional $C^*$-algebra $B_n$ and the norm closure of $\bigcup_n A_n$ is isomorphic to $A$. Furthermore we may assume, by modifying injections $A_n \subset A_{n+1}$, that for any natural number $k$ and $N$ there exists a non-zero projection $q \in A\cap B_k'$ such  that $ (\Ad u_k )^j(q)$, $j=0,1,...,N-1$ are mutually orthogonal projections, where $u_k$ is the canonical unitary of $1_{B_k}\otimes C(\T) \subset A_k$ (see \cite{Kis2}). 

For the finite subset $F \subset A^1$, we obtain $k$, $M\in \N$ and a finite subset $G$ of $(B_k)^1$ such that 
\[ F \subset_{4^{-1}\varepsilon} \{\sum_{n=-M}^{M} b_nu_k^n ; b_n \in G\}
\subset A_k ,\] where $A\subset_{\varepsilon}B$ means that for any $a\in A$ there exists $b \in B$ such that $\|a-b\|<\varepsilon$. Let $\delta = (2M(M+1))^{-1} \varepsilon $ and $N\in \N$ be such that $2N^{-1/2}< \delta$. Then, assumed as above, there exists a non-zero projection $q\in A\cap B_k'$ such that $(\Ad u_k )^j (q),\ j=0,1,...,N-1$ are mutually orthogonal projections in $A$.
It suffices to show that, for any $n\in \N$, there exist $m\geq n$ and $w\in B\cap B_k'$ satisfying the following conditions: $\|[w,u_k]\| < \delta $, and  
\[ \|w^*w- e_{m}\|<\varepsilon,\quad \|ww^*-ww^*f_m\|< \varepsilon,\quad \|w^2\| < \varepsilon. \] 
Indeed, for any $f\in F$ there exists $b_n \in G$ such that $\displaystyle 
\|f - \sum_{n=-M}^{M} b_nu_{k}^n \| < 4^{-1} \varepsilon$ and $\| b_n\| \leq 1$, then it follows that  
\begin{eqnarray}
\|[w,f]\| &<& 2^{-1}\varepsilon + \|w(\sum_{n=-M}^M  b_nu_{k}^n)-(\sum_{n=-M}^M  b_nu_{k}^n)w\| \nonumber \\
  &\leq &2^{-1}\varepsilon + \sum_{n=-M}^M \| [w, b_n u_k^n]\|< 2^{-1} \varepsilon + M(M+1)\delta = \varepsilon . \nonumber
\end{eqnarray}

Let $B_{k,q}$ denote the finite dimensional $C^*$-subalgebra $C^*(B_k \cup \{q\})$ of A. Since $(\Phi_{B_{k,q}} (e_n))_n = (e_n)_n$ and $ (\Phi_{B_{k,q}} (f_n))_n = (f_n)_n$, 
where $\Phi_{B_{k,q}}$ is the conditional expectation from $B$ into $B\cap B_{k,q}'$, we may assume that $e_n$, $f_n \in P(B \cap B_{k,q}')$ for any $n\in N$.
Since $A$ is a simple $C^*$-algebra, $\tau \in T(A) $ is faithful as a state. 
We denote by $\{e_{i,j}^{(l)}\}$ a system of matrix units of $B_k \cong \bigoplus_{l=1}^{L}M_{k_l}$. For any $l=1,2,...,L$ and $\tau \in T(A) $ we have $\tau (qe_{1,1}^{(l)})>0$. Since $\{\tau(qe_{1,1}^{(l)}); \tau \in T(A) \}$ is a compact subset of $(0,1]$ , we set $c_0 =\min \{ \tau (qe_{1,1}^{(l)}) ; \tau \in T(A), l=1,2,...,L\} >0 $. 

Set $\varepsilon_1 = \min\{(2N)^{-1}\varepsilon,\ 2^{-1}cc_0\}$. By the assumption for $(e_n)$ and $(f_n)$, there exists a large natural number $n_0$ satisfying the following conditions: $\|[ e_{n_0},x]\| < \varepsilon_1$, $\|[f_{n_0},x]\| <\varepsilon_1 $ for any $x\in\{e_{i,j}^{(l)}\}\cup\{u_k^j\}_{j=1,2,...,N-1}$, $c\tau(qe_{1,1}^{(l)})-\varepsilon_1 < \tau(f_{n_0}qe_{1,1}^{(l)})$ for any $\tau \in T(B) $ and any $l=1,2,...,L$, and $\tau (e_{n_0}) < 4^{-1} cc_0$ for any $\tau \in T(B)$. 
Then we have that for any $\tau \in T(B)$ and $l=1,2,...,L$, 
\[\tau(f_{n_0}qe_{1,1}^{(l)}) > c \tau (qe_{1,1}^{(l)}) -\varepsilon_1 \geq cc_0-\varepsilon_1 \geq 2^{-1}cc_0> 2\tau (e_{n_0}) \geq 2\tau(e_{n_0}e_{1,1}^{(l)}).\] 
Since $B$ is of real rank zero and the order on projections over $B$ is determined by traces, by Lemma \ref{lem:Orthogonal} there exists $w_{n_0}^{(l)} \in B$ such that 
\[ {w_{n_0}^{(l)}}^*w_{n_0}^{(l)} =e_{n_0} e_{1,1}^{(l)},\quad w_{n_0}^{(l)}{w_{n_0}^{(l)}}^* \leq f_{n_0}q e_{1,1}^{(l)},\quad \|(w_{n_0}^{(l)})^2\|< \varepsilon_1. \]
We define $w' \in B$ by
\[ w' = \sum_{l=1}^{L} \sum_{i=1}^{k_l} e_{i,1}^{(l)} w_{n_0}^{(l)} e_{1,i}^{(l)}.\]
Then we have that 
\[{w'}^*w' = e_{n_0},\quad w'{w'}^* \leq f_{n_0}q,\quad \| w'^2\| < \varepsilon_1,\]
and $[e_{i,j}^{(l)}, w']=0$; thus $w' \in B \cap B_k'$.
We define $w \in B$ by 
\[ w = N^{-1/2} \sum_{j=0}^{N-1} u_k^j w'u_k^{-j}.\]
Then we have that $w \in B \cap B_k' $ and  
\[ \|[w,u_k]\|=\| \Ad u_k (w) -w\| = N^{-1/2}\|(\Ad u_k)^N (w')-w'\|\leq 2N^{-1/2}<\delta.\] 
Since $(\Ad u_k)^j(q)$ are mutually orthogonal and $\|[e_{n_0},u_k^i]\|<\varepsilon_1$, we have that
\begin{eqnarray}
\|w^*w-e_{n_0}\|&=& \|N^{-1} \sum _{i=0}^{N-1} \sum_{j=0}^{N-1}u_k^i{w'}^*qu_k^{-i}u_k^jqw' u_k^{-j} -e_{n_0}\| \nonumber \\
&=& \| N^{-1} \sum _{i=0}^{N-1} u_k^i{w'}^*w'u_k^{-i} -e_{n_0}\| \nonumber \\
&\leq& \varepsilon_1< \varepsilon. \nonumber
\end{eqnarray}
Since $\| [ f_{n_0}, u_k^i] \| <\varepsilon_1$ and \[ww^* = N^{-1} \sum_{i,j} u_{k}^i w' u_{k}^{-i+j} {w'}^* f_{n_0} u_{k}^{-j},\] we have that $\|ww^* -  ww^*f_{n_0}\|$
\begin{eqnarray}
&\leq& N^{-1}\sum_{i,j}\|u_{k}^i  w' u_{k}^{-i+j}  {w'}^* f_{n_0} u_{k}^{-j} - u_{k}^i  w' u_{k}^{-i+j} {w'}^*  u_{k}^{-j} f_{n_0} \| \nonumber \\
&<& N\varepsilon_1 \leq \varepsilon. \nonumber 
\end{eqnarray}
Since $\|[e_{n_0}, u_k^i] \| < \varepsilon_1$ and $\|w'^2\| < \varepsilon_1$, we have that 
\[\|w^2\|\leq N^{-1} \sum_{i,j}\|u_k^{i}w'e_{n_0} u_k^{-i+j} w' u_k^{-j} \| < 2N\varepsilon_1 \leq \varepsilon.\]
\eprf

\begin{proposition}\label{prop:crossedproduct}
Let $A$ be a separable unital $C^*$-algebra with a unique tracial state and $\alpha$ be an automorphism of $A$ with the Rohlin property. 
Suppose that the order on projections over $A$ is determined by traces, $A$ has the property (SI), and the order on projections over $A \times_{\alpha} \Z$ is determined by traces. Then $A\times_{\alpha} \Z$ also has the property (SI) .
\end{proposition}
\bprf
Let $B$ be a unital $C^*$-algebra of real rank zero such that the order on projections over $B$ is determined by traces and $A\times_{\alpha}\Z \subset_{unital} B$. Let $(e_n)$, $(f_n) \in P(B^{\infty} \cap (A\times_{\alpha} \Z)')$ and  $c>0$ be as in the definition of (SI), i.e., they satisfy the following conditions: for any $p\in P(A\times_{\alpha}\Z) $ and any $\varepsilon>0$ there exists $N_p\in\N$ such that  $ c \tau(p)-\varepsilon \leq \tau(pf_np)$, $ n \geq N_p$, $ \tau \in T(B)$, and $\tau (e_n) \rightarrow 0 $ for any $\tau \in T(B)$.

Let $F$ be a finite subset of $(A\times_{\alpha}\Z)^1$ and $\varepsilon > 0 $. We obtain a finite subset $F_0$ of $A^1$ and $M \in \N$ such that 
\[ F \subset_{4^{-1}\varepsilon} \{ \sum_{n=-M}^M x_n u_{\alpha}^n : x_n \in F_0 \},\] where $u_{\alpha}$ is the canonical unitary of $A\times_{\alpha} \Z$ implementing $\alpha $.  Set $\delta = (2(M^2+3M+1))^{-1} \varepsilon $. It suffices to show that for any $n\in \N$ there exist $w \in B$ and $m\geq n$ satisfying the following conditions:  $\|[w,x]\| < \delta $ for any $x \in F_0 \cup \{ u_{\alpha}\}$, $\|w^*w- e_{m}\|<\varepsilon$, $ \|ww^*-ww^*f_m\|< \varepsilon$, and $ \|w^2\| < \varepsilon$. 
Indeed, then for any $f\in F$ there exists $x_n \in F_0$ such that $\displaystyle 
\|f - \sum_{n=-M}^{M} x_n u_{\alpha}^n \| < 4^{-1} \varepsilon$, and we have 
\[
\|[w,f]\| <2^{-1}\varepsilon + \sum_{n=-M}^M \| [w, x_nu_{\alpha}^n]\|< 2^{-1} \varepsilon + (M^2+3M+1)\delta = \varepsilon . \]
We let $N\in \N$ be such that $N > 4 \delta^{-2}$, $\tau_0$ be the unique tracial state of $A$, and $\displaystyle G=\bigcup_{k=0}^N \alpha^{-k} (F_0)$.
Since  $\alpha$ has the Rohlin property we can obtain a central sequence of projections $(p_{n})_n $ in $A$  such that $(\alpha ^j (p_{n}))_n$ are mutually orthogonal in $P(A_{\infty})$ for $\ j=0,1,...,N-1$ and $\tau_0 (p_n) \rightarrow 1/N$. Hence for any $q\in P(A)$ and $\varepsilon'>0$  there exists $n_0\in \N $ such that for any $ n\geq n_0$ and any $\tau \in T(B) $\[\tau (p_{n}) > N^{-1} - \varepsilon',\quad  |\tau(qp_{n})-\tau(q)\tau(p_{n})|< \varepsilon'. \]
Thus we have that for any $q\in P(A)$ and $\varepsilon'>0 $ there exists a large natural number $n_q$ such that for any $n \geq n_q$ and any $\tau \in T(B) $, 
\[\tau (q p_{n} q) \geq N^{-1}\tau(q)- \varepsilon'.\]
By the condition (1) for $f_n$ and the fact $(f_n)\in B^{\infty} \cap (A\times_{\alpha} \Z)' $, we inductively obtain $l_n \in \N$ such that $l_n >l_{n-1}$, $\tau(p_nf_{l_n}p_n) \geq c \tau (p_n)- n^{-1}$ for any $\tau \in T(B)$, and $\|[p_n, f_{l_n}]\| < n^{-1}$.  
We define a projection $\widetilde{f}_n$ from $p_nf_{l_n}p_n$ by continuous function calculus such that $\widetilde{f}_n \leq p_n$ and $\| \widetilde{f}_n -p_n f_{l_n}p_n \| \rightarrow 0$ and set $\widetilde{e}_n= e_{l_n}$. Note that $(\widetilde{f}_n) \in P(B^{\infty}\cap A')$.  
Since $A$ has a unique tracial state, we have that for any $q \in P(A) $ and $\tau \in T(B) $, 
\[\lim_{n\rightarrow \infty}\tau(\widetilde{f}_n)^{-1}\tau(\widetilde{f}_n q \widetilde{f}_n) = \tau(q) ,\] 
which implies that for any $q \in P(A)$ and $ \tau \in T(B) $,
\[\lim_{n \rightarrow \infty} |\tau(q \widetilde{f}_nq) - \tau(q) \tau(\widetilde{f}_n)|=0.\]  
Thus for any $q\in A$ and any $\varepsilon' >0$ there exists $N_q \in \N$ such that for $n\geq N_q$ and $\tau \in T(B)$,
\begin{eqnarray}
\tau(q\widetilde{f}_nq)& >& \tau(\widetilde{f}_n)\tau(q)-\varepsilon'> (c\tau(p_n)-\varepsilon')\tau(q) -\varepsilon'\nonumber\\
&>& cN^{-1}\tau (q) - \varepsilon'((c+1)\tau(q)+1). \nonumber
\end{eqnarray} 
Then we can apply the property (SI) of $A$ to $(\widetilde{e}_n)$, $(\widetilde{f}_n)\in P(B^{\infty}\cap A')$, and for $\varepsilon_1=(3N)^{-1}\varepsilon$ we  obtain $w' \in B$ and $m \geq n$
such that $m^{-1}< \varepsilon $, $\|[w', g]\| < N^{-1/2}\delta$ for any $g \in G \subset A$,
\[ \|(w')^*w' - \widetilde{e}_{m}\|<\varepsilon_1,\quad \|w'(w')^* -w'(w')^*\widetilde{f}_m\|< \varepsilon_1,\quad \|(w')^2\|<\varepsilon_1 .\]
Moreover, by $(e_n)$, $(f_n) \in B^{\infty}\cap (A\times _{\alpha} \Z)' $, by the property on $(p_n)$, and by $\| [ p_n, f_{l_n}]\| < n^{-1},$ we may assume that for any $j=1,2,...,N$,
\[ \|[\widetilde {e}_m, u_{\alpha}^j]\| < \varepsilon_1,\quad \|[f_{l_m}, u_{\alpha}^j]\| < \varepsilon_1,\quad \|\alpha^j(p_m)p_m\| < \varepsilon_1,\quad \|\widetilde{f}_m - p_m f_{l_m} \| < \varepsilon_1.\]  By replacing $w'$ with $\widetilde{f}_m w' \widetilde{e}_m $ we may further assume that 
\[w'{w'}^* \leq \widetilde{f}_{m},\quad {w'}^*w' \leq \widetilde{e}_m.\]
We define $w$ as $\displaystyle N^{-1/2} \sum_{j=0}^{N-1} (\Ad u_{\alpha} )^j (w') \in B $. Then we have 
\[\|[w, u_{\alpha}]\|=\|\Ad u_{\alpha} (w)-w\|= N^{-1/2}\| \Ad u_{\alpha}^N (w')-w'\| <2N^{-1/2} <\delta,\] 
and for any $f \in F_0$ ,
\begin{eqnarray} 
\|[w, f] \| &=& N^{-1/2}\|\sum_{j=0}^{N-1}[\Ad u_{\alpha}^j(w'),f] \| \nonumber\\
&\leq& N^{-1/2} \sum \|[w', \Ad u_{\alpha}^{-j}(f)] \| <\delta. \nonumber
\end{eqnarray} 
Since $w'{w'}^* \leq \widetilde{f}_m \leq p_m$, $\|\alpha^j(p_m)p_m\| <\varepsilon_1$, $\|{w'}^*w' - \widetilde{e}_{m}\|<\varepsilon_1$, and $\|[\widetilde {e}_m, u_{\alpha}^j]\| < \varepsilon_1$ we have that 
\begin{eqnarray}
\|w^*w - e_{l_m}\|&=&\|N^{-1}\sum_{i,j}u_{\alpha}^i w'^* p_m u_{\alpha}^{-i+j} p_m w'u_{\alpha}^{-j} -\widetilde{e}_m\| \nonumber\\
&<& \| N^{-1}\sum_i u_{\alpha}^i {w'}^*w' u_{\alpha}^{-i} -\widetilde{e}_m \|
+ (N-1)\varepsilon_1 \nonumber \\
& <& \| N^{-1}\sum_i u_{\alpha}^i \widetilde{e}_m u_{\alpha}^{-i} -\widetilde{e}_m \| + N\varepsilon_1 \nonumber \\
&<&(N+1) \varepsilon_1< \varepsilon .\nonumber
\end{eqnarray}
Since $\|\widetilde{f}_m - p_m f_{l_m} \| < \varepsilon_1$, we have that 
\[\|(w')^* - (w')^*f_{l_m}\| = \| (w')^* \widetilde{f}_m - (w')^* p_m f_{l_m} \| < \varepsilon_1 ,\] and since $\|[f_{l_m}, u_{\alpha}^j]\| < \varepsilon_1$ we have that $\|ww^* -  ww^*f_{l_m}\|$
\begin{eqnarray}
&\leq& N^{-1}\sum_{i,j}\|u_{\alpha}^i  w' u_{\alpha}^{-i+j}  {w'}^*  u_{\alpha}^{-j} - u_{\alpha}^i  w' u_{\alpha}^{-i+j} {w'}^*  u_{\alpha}^{-j} f_{l_m} \| \nonumber \\
&<& N^{-1}\sum_{i,j}\|u_{\alpha}^i  w' u_{\alpha}^{-i+j}  {w'}^* f_{l_m} u_{\alpha}^{-j} - u_{\alpha}^i  w' u_{\alpha}^{-i+j} {w'}^*  u_{\alpha}^{-j} f_{l_m} \| + N \varepsilon_1 \nonumber \\
&<& 2N\varepsilon_1 < \varepsilon. \nonumber 
\end{eqnarray} 
Since $\|[\widetilde{e}_m, u_{\alpha}^j]\| < \varepsilon_1$, $\| (w')^*w'-\widetilde{e}_m\| < \varepsilon_1 $ and $\| w'^2 \| < \varepsilon_1$, we have that 
\begin{eqnarray}
\|w^2\| &\leq& N^{-1} \sum_{i,j} \| u_{\alpha}^iw'\widetilde{e}_m u_{\alpha}^{-i+j}w' u_{\alpha}^{-j}\| \nonumber \\ 
&<& N^{-1} \sum \|u_{\alpha}^iw'u_{\alpha}^{-i+j}\widetilde{e}_m w' u_{\alpha}^{-j}\| + \varepsilon_1 N  \nonumber \\ 
&<& N^{-1} \sum \|u_{\alpha}^iw'u_{\alpha}^{-i+j} (w')^*w' w' u_{\alpha}^{-j}\| + 2N\varepsilon_1   < 3N \varepsilon_1 = \varepsilon. \nonumber
\end{eqnarray}  \eprf

On the assumption for the crossed product in the above proposition, Osaka and Phillips showed in \cite{OP1} that: if $A$ is an infinite dimensional stably finite simple unital $C^*$-algebra of real rank zero, the order on projections over $A$ is determined by traces, and $\alpha \in \Aut (A) $ has the tracial Rohlin property, then the order on projections over $A\times_{\alpha} \Z $ is determined by traces.

In \cite{OP2}, Osaka and Phillips proved that the crossed product of a unital simple TAF-algebra with a unique tracial state by an automorphism of $A$ with the
tracial Rohlin property also has a unique tracial state.
In \cite{Lin1}, H.Lin proved that the crossed product of a unital simple TAF-algebra with a unique tracial state by an automorphism of $A$ with the cyclic Rohlin property is also a TAF-algebra. 
If $A$ is a TAF-algebra, it follows from \cite{Lin2} that the order on projections over $A$ is determined by traces.
Using these results, we can see the following corollary.
\begin{corollary}\label{cor:cross}
Let $A$ be a separable unital simple $C^*$-algebra with a unique tracial state and $\alpha$ an automorphism of $A$ with the Rohlin property.
If A is a TAF-algebra having the property (SI),
then $A\times_{\alpha}\Z$ is a TAF-algebra with a unique tracial state having the property (SI).
\end{corollary}

Before closing this section, we show that the class of TAF-algebras with a unique tracial state having the property (SI) is closed under inductive limits.
\begin{proposition}\label{prop:inductivelimit}
Let ($A_n$, $\varphi_n$) be an inductive sequence of separable unital simple $C^*$-algebras, where $\varphi_n$ is a unital $*$-homomorphism from $A_n$ into $A_{n+1}$. If for any $n\in\N$, the order on projections over $A_n$ is determined by traces
and $A_n$ has the property (SI), then the inductive limit $C^*$-algebra $\displaystyle \lim_{\longrightarrow} (A_n, \varphi_n)$ also has the property (SI). 
\end{proposition}
\bprf
 Let $A$ be the inductive limit $C^*$-algebra of the system $(A_n, \varphi_n)$, and $\widetilde{\varphi}_n$ the canonical map from $A_n$ into $A$.
Assume that $A_n$ has the property (SI) for any $n\in \N$ and assume that a $C^*$-algebra $B$, $(e_n)_n$, and $(f_n)_n \in B^{\infty}\cap A'$ satisfy the assumptions of the property (SI).
Let $F_n$, $n\in \N$ be a finite subsets of $A^1$ such that $F_n \subset {\widetilde{\varphi}_n (A_n)}^1$ and $\bigcup F_n$ is a dense subset of $A^1$.
By the property (SI) of $A_n$, we inductively obtain $w_n \in B$, and $m_n\in \N$ such that $\|{w_n}^*w_n-e_{m_n}\| < n^{-1}$, $\|w_nw_n^*- w_nw_n^*f_{m_n} \| < n^{-1}$, $\|{w_n}^2\|<n^{\-1}$, $\|[w_n,a]\|<n^{-1}$ for any $a\in F_n $, and $m_n >m_{n-1} $.
Then we have that $(w_n)_n\in B^{\infty} \cap A'$, $(w_n)^*(w_n)=e_{m_n}$, $(w_n)(w_n)^*\leq f_{m_n}$, and $(w_n)^2=0$. \eprf

The inductive limit $C^*$-algebra of a sequence of unital simple TAF-algebras with a unique tracial state is also unital simple TAF-algebra with a unique tracial state. Thus we obtain the following corollary.
\begin{corollary}\label{cor:inductivelimit}
Let ($A_n$, $\varphi_n$) be an inductive sequence of unital $C^*$-algebras with unital homomorphisms $\varphi_n$ and let $A$ be the inductive limit $C^*$-algebras of ($A_n$, $\varphi_n$). If $A_n$ is a unital simple TAF algebra with a unique tracial state and has the property (SI) for any $n\in \N$, then so is $A$.

\end{corollary}

\section{Realizing automorphisms with the Rohlin property}\label{sec:M}
In this section, using a certain automorphism of the Jiang-Su algebra which is aperiodic in $\Aut(\Zj)/\{\Ad W\in \Aut(\Zj): W\in U(\pi_{\tau}(\Zj)'')\}$ and using the property (SI), we shall obtain automorphisms of certain $C^*$-algebras with the Rohlin property.

Let $\Zj$ be the Jiang-Su algebra defined in \cite{JS}.
We first remark that $\Zj$ is isomorphic to $\bigotimes_{n\in \N}\bigotimes_{i=1}^{n} \Zj $ (\cite{JS}) and we define an automorphism $\sigma_n$ of $\displaystyle \Zj _{n}=\bigotimes_{i=1}^n \Zj$ by 
\[\sigma_n( x_1\otimes x_2 \otimes \cdots \otimes x_N )=x_2\otimes x_3\otimes \cdots \otimes x_{N}\otimes x_1.\]
We then define an automorphism $\sigma $ of $\Zj=\bigotimes_{n\in \N} \Zj_n$ by $\sigma=\bigotimes_{n\in N} \sigma_n$. It follows that $\sigma$ is an automorphism of $\Zj$ whose non-zero powers are not weakly inner in $\pi_{\tau}$.

A similar argument in the next lemma appeared in the proof of Lemmas 6.3, 6.4, 6.5, and 6.6 of \cite{RW}. 

\begin{lemma}\label{lem:measure}
Let $\mu$ be a probability measure on the $\sigma$-algebra of Borel sets of $[0,1]$ such that $f \mapsto  \int f d\mu $ is a faithful state on C([0,1]), $\mu_n$ be the product measure $\mu \times \mu \times \cdots \times \mu $ on $[0,1]^n$, and $\sigma_n$ be the automorphism of $C([0,1]^n)$ defined by 
\[\sigma_n(x_1 \otimes x_2 \otimes \cdots \otimes x_n) = x_2 \otimes x_3 \otimes \cdots \otimes x_n\otimes x_1,\quad x_i \in C([0,1]) .\]
For any $k \in \N$ and $\varepsilon >0 $, there exists $ m \in \N$ and mutually orthogonal positive elements $f_0,f_1,...,f_{k-1} \in C([0,1]^m)_+^1$ such that 
\[ \sigma_m (f_j) =f_{j+1}, \quad j=0,1,...,k-2,\quad{\rm and}\quad \int_{[0,1]^m} (1-\sum_{j=0}^{k-1} f_j) d \mu_m < \varepsilon. \]
\end{lemma}   
\bprf
For any probability measure $\mu $ on $[0,1]$ as in the assumption, we obtain $g_{n}^{(0)},g_{n}^{(1)} \in C([0,1])^1_+$, $n\in \N$ and a constant $c_{\mu}<1$ such that 
\[g_{n}^{(0)}g_{n}^{(1)}=0,\ n\in \N,\quad \lim_{n\rightarrow \infty} \mu (g_{n}^{(0)}+g_{n}^{(1)})=1,\]\[ c_{\mu} \geq \max \{\mu(g_{n}^{(l)}): l=0,1,\quad n\in\N\}, \]
where $\mu (f) $ means that $\int f d\mu$.
Remark that $g_n^{(0)}$ and $g_n^{(1)}$ are non-zero elements.
Set $\displaystyle d_n = g_{n}^{(0)}+ g_{n}^{(1)}$ and let $m$ be a large prime number such that $2k/m < \varepsilon $ and $(c_{\mu})^m < \varepsilon/8$ and let $n_0 \in \N$ be such that $(\mu(d_{n_0}))^m > 1- \varepsilon /4.$

We define   
\[S=\{ (l_1,l_2,...,l_m) : l_i \in \{ 0,1\},\quad \exists i,i' \in\{1,...,m \}\ {\rm such\ that\ } l_i \neq l_{i'} \}. \]
Since $m$ is a prime number, $x=g_{n_0}^{(l_1)}\otimes g_{n_0}^{(l_2)} \otimes \cdots \otimes g_{n_0}^{(l_m)}$ with $(l_1,l_2,...,l_m) \in S$ satisfies that $\sigma_m^i(x)x=0$ for any $i=1,2,...,m-1$. We define an equivalent relation $\sim$ on $S$ by $(l_1,l_2,...,l_m) \sim (l_1',l_2',...,l_m')$ if there is $i\in\{1,2,...,m-1\}$ such that $\sigma_m^i(g_{n_0}^{(l_1)},\otimes \cdots \otimes g_{n_0}^{(l_m)})=g_{n_0}^{(l_1')}\otimes \cdots \otimes g_{n_0}^{(l_m')}$, and define positive elements $h_i\in C([0,1]^m)$ by 
\[h_1= \sum_{(l_1,...,l_m) \in S/ \sim } g_{n_0}^{(l_1)}\otimes g_{n_0}^{(l_2)}\otimes \cdots \otimes g_{n_0}^{(l_m)},\quad h_i= \sigma^{i-1}(h_1),\]
where the summation is taken by choosing a representative for each equivalence class of $S/\sim$.
Remark that $h_i$, $i=1,2,...,m-1$ are mutually orthogonal and $\mu_m (h_i)=\mu_m(h_1)$.
Since \[d_{n_0}\otimes\cdots \otimes d_{n_0} = \sum_{i=1}^{m}h_i + \sum_{j=0,1} g_{n_0}^{(j)}\otimes g_{n_0}^{(j)}\otimes \cdots \otimes g_{n_0}^{(j)},\]
we have that 
\begin{eqnarray}
\mu_{m}(\sum_{i=1}^m h_i) &=& \mu (d_{n_0})^m - \sum_{j=0,1}\mu(g_{n_0}^{(j)})^m \nonumber \\
&\geq& \mu(d_{n_0})^m-2(c_{\mu})^m  \nonumber \\
&>& \mu (d_{n_0})^m - \varepsilon /4 = 1- \varepsilon /2 .\nonumber 
\end{eqnarray}

Let $a$ and $r$ be the natural numbers such that $m=ak+r$ and $0\leq r < k$ and define $f_j \in C([0,1]^m)^1_+$, $j=0,...,k-1$ by 
\[f_0=h_1+h_{1+k}+h_{1+2k}+\cdots + h_{1+(a-1)k}, \quad f_j=\sigma_m^j(f_0).\]
Thus we have that $f_if_j=0$, $i \neq j$ and 
\begin{eqnarray}
\mu_m(\sum_{j=1}^k f_j) &=& \mu_m (\sum_{i=1}^{ak} h_i) = \mu_m (\sum_{i=1}^m h_i) -\mu_m(\sum_{i=ak+1}^m h_i)\nonumber \\
 &>& 1-\varepsilon/2  -r\mu_m(h_1)>1-\varepsilon/2-k/m > 1- \varepsilon. \nonumber
\end{eqnarray}
\eprf

\begin{lemma}\label{lem:p2}
Let $A$ be a unital TAF-algebra with $A\otimes\Zj \cong A$. Then for any automorphism 
$\alpha $ of $A$ and any $k\in \N$, there exists a $(p_n)_n \in P((A\otimes \Zj) _{\infty}) $ satisfying the following conditions : $((\alpha
\otimes \sigma)^j(p_n))_n$ are mutually orthogonal as projections in $(A\otimes \Zj )_{\infty}$ for $j=0,1,...,k-1$, 
\[\tau(1_{A\otimes \Zj}-\sum_{j=0}^{k-1} (\alpha
\otimes \sigma)^j(p_n))<n^{-1}, \quad \tau \in T(A\otimes \Zj), \]
 and there exists $c>0$ satisfying that for any $q\in P(A\otimes \Zj )$ and $\varepsilon > 0$ there exists $N_q \in \N$ such that 
\[c\tau (q)- \varepsilon \leq \tau(qp_n),\quad \tau \in T(A\otimes \Zj),\ n\geq N_q .\]
\end{lemma}
\bprf
Let $F$ be a finite subset of $(A\otimes \Zj)^1$, $G$ a finite subset of $P(A\otimes\Zj)$, $\varepsilon>0$, $k\geq 2$, and $\varepsilon_1= \frac{\varepsilon}{6k} $.
Since $A \otimes \Zj$ is a TAF-algebra, there exists a finite dimensional $C^*$-subalgebra $B$ of $ A \otimes \Zj$ satisfying that $\|[1_B,x]\| < \varepsilon_1 $ for any $x\in F\cup G$,
$\{x 1_B : x \in F \cup G\} \subset _{\varepsilon_1} B^1$, and 
\[\tau (r)< \varepsilon_1 ,\quad \tau \in T(A\otimes \Zj), \]
where $r=1_{A\otimes \Zj} -1_B$. 
Let $\{ e_{i,j}^{(l)}\}$ be a system of matrix units of $B \cong \bigoplus_{l=1}^{L} M_{k_l} $. Let $\delta=(16k\dim(B))^{-1}\varepsilon $ and let $N$ be such that $\displaystyle \{e_{i,j}^{(l)} \} \subset _{\delta} A\otimes \bigotimes_{n=1}^{N} \Zj_n \subset
A\otimes \Zj $.  
It suffices to show that there exists a projection $p_0 \in (1-r)A\otimes \Zj (1-r)\cap B'$ such that $\|(\alpha \otimes \sigma)^j (p_0)p_0 \|< \varepsilon$ for $j=1,2,...,k-1$,  
\[ \tau(1_{A\otimes \Zj} -\sum_{j=0}^{k-1} (\alpha \otimes \sigma)^j (p_0)) <\varepsilon,\quad k^{-1}\tau(q)- \varepsilon < \tau(qp_0)\] for any $\tau \in T(A\otimes \Zj)$, and any $q\in G$.
This is shown as follows. First, we have that for any $f \in F$
\begin{eqnarray}
\|[p_0, f] \| &\leq &\| p_0 f(1-r)- f(1-r)p_0 \| + \varepsilon_1 \nonumber \\
&\leq & \| p_0 b - bp_0 \| +  3\varepsilon_1, \nonumber
\end{eqnarray}
for some $b \in B^1$. Thus $\| [ p_0 , f ] \| < \varepsilon $.
Then for an increasing sequence $(F_n)$ of finite subsets of $(A\otimes \Zj) ^1$
with $ \overline{\bigcup F_n}= (A \otimes \Zj) ^1 $ and an increasing sequence $(G_n)$ of finite subsets of $P(A \otimes \Zj ) $ with $ \overline{\bigcup G_n} = P(A \otimes \Zj) $, we can inductively  obtain $p_n \in P(A\otimes \Zj) $ such that $\|[p_n,f]\| <n^{-1}$ for any $f \in F_n$, $\|(\alpha\otimes\sigma)^j(p_n)
p_n \| < n^{-1}$ for any $j=1,2,...,k-1$,
$ \tau(1_{A\otimes \Zj } - \sum_{j} (\alpha\otimes \sigma)^{j} (p_n)) < n^{-1}$ for any $\tau \in T(A\otimes \Zj) $, and $k^{-1} \tau(q) - n^{-1} < \tau(qp_n)$ for any $\tau \in T(A\otimes \Zj)$ and any $q\in G_n$. 
Thus $(p_n)_n$ would satisfy the desired conditions in the lemma.

Because $C([0,1])$ is embedded in $\Zj$ unitally, by Lemma \ref{lem:measure}, there is a natural number $m$ with $m\geq \max \{k,N\} $ and mutually orthogonal positive elements 
$f_0, f_1, ..., f_{k-1} \in {\Zj_m}^1$ such that 
\[ \sigma_m(f_j)=f_{j+1}, \quad j=0,1,...,k-2,\]
\[\tau_1 (1- \sum_{j=0}^{k-1} f_j)< \varepsilon/4 , \]
where $\tau_1$ is the unique tracial state on $\Zj_m$. Remark that $\tau_1 (f_j)= \tau_1(f_0) > k^{-1}(1-\varepsilon/4)$.
We define $\widetilde{f}_j$ as the image of $f_j$ under the natural embedding $
\Zj_m \rightarrow 1_A\otimes \bigotimes_{n\in\N} \Zj_n $. Set $D=A\otimes \Zj \cap B'$ and $\widetilde{g}_0 = \Phi_{B+\C r}(\widetilde{f}_0) \in D_+$, and 
\begin{eqnarray}
f_{\varepsilon_1}(t)=\left\{ \begin{array}{ll}
t/{\varepsilon_1},\quad & 0\leq t < \varepsilon_1 \\
1,\quad  & \varepsilon_1 \leq t< 1. \\
\end{array} \right.\nonumber
\end{eqnarray}   
If we choose sufficiently large $m$, then $\| \widetilde{g}_0-\widetilde{f}_0\| $ goes to $0$. Thus we may assume that 
\[ \| \widetilde{g}_0 - \widetilde{f}_0 \| < \delta \quad {\rm and} \quad \|f_{\varepsilon_1}(\widetilde{g}_0)-f_{\varepsilon_1}(\widetilde{f}_0) \| < \varepsilon/4.\]

Since $D$ has real rank zero, $D_{r,0}:=\overline{(\widetilde{g}_0 D(1-r) \widetilde{g}_0)}$ also has real rank zero. Thus there exists a positive element $g_{0} \in D_{r,0}$ with finite spectrum such that 
\[g_{0}=\sum_{i=1}^{M}\lambda_{0,i} p_{0,i},\ p_{0,i}\in P(D_{r,0}),\ 0 < \lambda_{0,i} \leq 1,\]
 \[ \|g_0 -\widetilde{g}_0(1-r)\|<\delta ,\  {\rm and}\  
 \| f_{\varepsilon_1} (g_0)- f_{\varepsilon_1} ((1-r)\widetilde{g}_0) \| <\varepsilon/4.\] 
We define a projection $p_0$ of $D_{r,0} \subset A \otimes \Zj \cap B' $ as $\sum_{\lambda_{0,i}\geq \varepsilon_1} p_{0,i} $. 
From $p_0= f_{\varepsilon_1} (g_0) p_0$, it follows that 
\[(\alpha \otimes \sigma)^j (p_0) p_0 = (\alpha \otimes \sigma)^j(p_0 )(\alpha \otimes \sigma)^j(f_{\varepsilon_1}(g)) f_{\varepsilon_1} (g) p_0.\]
Since $\|f_{\varepsilon_1} (g) - f_{\varepsilon_1}((1-r)\widetilde{g}_0)\| < \varepsilon /4 $ and $\| f_{\varepsilon_1} (\widetilde{g}_0) - f_{\varepsilon_1}(\widetilde{f}_0) \| < \varepsilon / 4$, we have that 
\begin{eqnarray}
 \| (\alpha\otimes \sigma)^j (f_{\varepsilon_1} (g)) f_{\varepsilon_1} (g) \|&\leq &  \| (\alpha \otimes \sigma)^j ((1-r)f_{\varepsilon_1}(\widetilde{g}_0)) f_{\varepsilon_1}(\widetilde{g}_0)(1-r) \| + \varepsilon/2 \nonumber \\
& \leq & \| (\alpha \otimes \sigma)^j (f_{\varepsilon_1}(\widetilde{f}_0)) f_{\varepsilon_1}(\widetilde{f}_0) \| + \varepsilon = \varepsilon .\nonumber
\end{eqnarray}
Thus $\|( \alpha\otimes\sigma )^j (p_0) p_0 \| < \varepsilon $.

Since $g_0 \leq p_0 + \varepsilon_1 1_{D_{r,0}}$, $\| g_0- \widetilde{g}_0(1-r) \|< \delta$, $\tau(r) < \varepsilon_1 $ for any $\tau \in T(A\otimes \Zj)$, $\|\widetilde{f}_0 -\widetilde{g}_0 \| < \delta$, and $\displaystyle \tau_1 (1-\sum_j f_j)< \varepsilon/4  $, we have that for any tracial state $\tau$ of $A\otimes \Zj $,
\begin{eqnarray}
\tau (1_{A\otimes \Zj} - \sum_{j=0}^{k-1} (\alpha \otimes \sigma)^j(p_0) )&\leq&\tau (1- \sum_j (\alpha \otimes \sigma )^j (g_0 ))+ k \varepsilon_1 \nonumber \\
&<&\tau (1-\sum_j(\alpha\otimes\sigma)^j(\widetilde{g_0}(1-r) ))+k (\delta +
\varepsilon_1) \nonumber \\
&<&\tau (1-\sum \widetilde{f}_j) +2k(\delta + \varepsilon_1) \nonumber  \\
& < &\tau_1(1_{\Zj_m}-\sum_j f_j )+ 2/3\varepsilon < \varepsilon . \nonumber
\end{eqnarray}
Let $q\in G$ and $\tau \in T(A\otimes \Zj)$. Since $G(1-r) \subset_{\varepsilon_1} B^1$, there is $b \in B^1 $ such that $\|b - q(1-r)\| < \varepsilon_1 $ and since $\{ e_{i,j}^{(l)}\} \subset_{\delta} A \otimes \bigotimes_{n=1}^N \Zj_n$, there is $a \in (A \otimes \bigotimes_{n=1}^N \Zj_n)^1$ such that $\|b-a\| < 2\dim (B) \delta$. 
Since $g_0 \leq p_0+ \varepsilon_1 1$, $\|g_0 - \widetilde{g}_0 (1-r) \|< \delta$, and $\|\widetilde{f}_0-\widetilde{g}_0\| < \delta$, we have that for any $q \in G$ and $\tau \in T(A\otimes \Zj) $ 
\begin{eqnarray}
\tau(qp_0) &\geq& \tau(q g_0) - \varepsilon_1 \geq  \tau(q(1-r)\widetilde{g}_0 (1-r) q ) -\delta-\varepsilon_1 \nonumber \\
&\geq & \tau(b \widetilde{g}_0 b^*) -(\delta + 3\varepsilon_1) \nonumber \\
&\geq & \tau(b \widetilde{f}_0 b^*) -(2\delta + 3\varepsilon_1) \nonumber \\
&\geq & \tau(a\widetilde{f}_0 a^*)-(4\dim(B)\delta + 2\delta + 3\varepsilon_1) \nonumber \\
&= & \tau_1(f_0)\tau(aa^*) - (4\dim(B)\delta + 2\delta + 3\varepsilon_1), \nonumber
\end{eqnarray}
where the last equality follows from $a\widetilde{f}_0a^* =aa^* \otimes f_0 \in (A\otimes \bigotimes_{n=1}^N \Zj_n) \otimes \Zj_m$. Since 
\begin{eqnarray}
|\tau(aa^*)-\tau(q)| &<& |\tau(aa^*)-\tau(bb^*)|+|\tau(bb^*)-\tau(q(1-r)q)|+|\tau(r)| \nonumber \\
&<& 4\dim(B)\delta + 3\varepsilon_1, \nonumber
\end{eqnarray}
we have that $\tau_1(f_0)\tau(aa^*) - (4\dim(B)\delta + 2\delta + 3\varepsilon_1) $
\begin{eqnarray}
&>& \tau_1(f_0)\tau(q) - (8\dim(B)\delta+2\delta+ 6\varepsilon_1) \nonumber \\
&> & k^{-1} \tau (q) - (4k)^{-1}\varepsilon - (8\dim(B)\delta+2\delta+ 6\varepsilon_1) \nonumber \\
&> & k^{-1}\tau(q) - \varepsilon/8- ( 1/4 + 1/16+ 1/2)\varepsilon > k^{-1}\tau(q)- \varepsilon . \nonumber 
\end{eqnarray}
Hence $\tau(qp_0)>k^{-1}\tau(q) - \varepsilon$.
\eprf

The basic idea of the proof of the following lemma was appeared in the proof of Lemma 4.6 in \cite{Kis2} for a simple unital A$\T$-algebra case. 

\begin{lemma}\label{lem:p3}
Let $A$ be a separable unital $C^*$-algebra of real rank zero and let $\alpha $ be an automorphism of $A$.
Suppose the order on projections over $A$ is determined by traces, $A$ has the property (SI),
and there exists a $(p_n) \in P(A_{\infty})$ satisfying the following conditions: all $p_n$ are projections,
\[(\alpha^i (p_n)) (\alpha^j (p_n))=0,\ i\neq j,\quad i,j=0,1,...,k-1 \]
\[\tau(1-\sum_{j=0}^{k-1} \alpha^j(p_n)) \rightarrow 0, \quad \tau \in T(A),\]
 and there exists a $c>0$ satisfying that for any $q\in P(A )$ and $\varepsilon > 0$ there exists $N_q \in \N$ such that 
\[c\tau (q)- \varepsilon \leq \tau(qp_nq),\quad \tau \in T(A),\ n\geq N_q .\]
Then for any finite subset $F$ of $A$, $\varepsilon>0$, and $k\in\N$,
there exist $e_0^{(0)}$, $e_1^{(0)}$,$...$,$e_{k-1}^{(0)}$,$e_0^{(1)}$,$e_1^{(1)}$,$...$,$e_k^{(1)}\in P(A) $ and $w\in U(C([0,1])\otimes A)$ such that $\sum e_i^{(0)} + \sum e_j^{(1)} =1_A$,
\[ w(0)=1_A,\quad \|[w,1_{C([0,1])}\otimes x]\|< \varepsilon,\ x\in F,\]
\[\| \Ad w(1) \circ \alpha (e_j^{(i)}) - e_{j+1}^{(i)}\|<\varepsilon ,\quad i=0,1, \quad j=0,1,...,k+i-2,\]
\[ \| [e_j^{(i)},x]\|<\varepsilon,\ i=0,1,\ j=0,1,...,k+i-1,\ x\in F.\]
\end{lemma}
\bprf
Let $k$ be a natural number and let $(p_{n})_n \in P(A_{\infty})$ be as in the lemma, and let $p_{j,n}=\alpha^j(p_n)$, $j=0,1,...,k-1$. We define $e'_n$ as $1_A- \sum_{j} p_{j,n}$.

 Since $(e_n')$ and $(\alpha(p_{k-1,n}))$, as elements of $P(B^{\infty} \cap A')$ with $B=A$, satisfy the conditions (1), (2) for $(e_n)$, $(f_n)$ in the property (SI),  there exists a $(w_n) \in A_{\infty}$ such that 
\[(w_n)^*(w_n)=(e'_{n}),\quad (w_n)(w_n)^* \leq \alpha ((p_{k-1,n})),\quad (w_n)^2=0 .\]
We define a $(u_n)_n \in U(A_{\infty})$ by $u_n =w_n+ w_n^* +
1 - w_n^*w_n - w_nw_n^* $, and  $(E_n)_n\in P(A_{\infty})$ 
by $E_n=e'_{n} + p_{0,n}$. Note that $(u_n)(1_{A_{\infty}}-(E_n))=1_{A_{\infty}}-(E_n)$ and $  (e'_n)_n \leq \Ad (u_n)_n \circ \alpha ((p_{k-1,n})_n) $. Since $(u_n)$ is a self adjoint element there exists 
a $(F_n)_n \in P(A_{\infty})$ such that $(u_n) = 2(F_n)-1 .$

We define  central sequences $(e_{j,n}^{(1)}) $, $j=0,1,...,k$ by $e_{j,n}^{(1)}=(\Ad u_n\circ \alpha)^{j-k} (e'_n)$. By
$(u_n)(1_{A_{\infty}}-(E_n))=1_{A_{\infty}}-(E_n)$, we have that $(p_{j,n})\geq(e_{j,n}^{(1)}),\ j=0,1,...,k-1$, and thus $(e_{i,n}^{(1)})(e_{j,n}^{(1)})=0,$ $i$\noteq$j$, $i,j=0,1,...,k$.
Let $e_{j,n}^{(0)}=p_{j,n}-e_{j,n}^{(1)}$ for $f=0,1,...,k-1$. Then we have $(e_{j,n}^{(0)}) \in P(A_{\infty}) $ and  
\[\Ad (u_n)\circ \alpha ((e^{(0)}_{j,n}))= (e^{(0)}_{j+1,n}),\quad j=0,1,...,k-2,\]
\[\Ad (u_n)\circ \alpha ((e^{(1)}_{j,n}))= (e^{(1)}_{j+1,n}),\quad j=0,1,...,k-1,\]
\[\sum_{j=0}^{k-1}(e^{(0)}_{j,n}) + \sum_{j=0}^{k}(e^{(1)}_{i,n})= 1_{A_{\infty}}.\]
By functional calculus we may assume that $e^{(0)}_{j,n}$, $e^{(1)}_{j,n}$ are mutually orthogonal projections of $A$ such that $\sum_j e^{(0)}_{j,n}+\sum_j e^{(1)}_{j,n}=1_A$ and
$F_n$ is a projection of $A$. We define a path of unitaries $w_n(t)$ by $w_n(t)=F_n+ \exp (\pi i t )(1-F_n ) $. If $n$ is sufficiently large then $w_n(t)$, $e^{(i)}_{j,n}$ satisfy the desired conditions.
\eprf

\begin{definition}\label{def:asym}
Let $A$ and $B$ be unital $C^*$-algebras. We say that unital $*$-homomorphisms $\varphi$ and $\psi$ from $A$ to $B$ are asymptotically unitarily equivalent if there is a continuous path of unitaries $u(t) \in U(B)$, $t\in[0,\infty)$ such that $\varphi (a)=\lim_{t \rightarrow \infty} \Ad u(t)\psi(a) $ for all $a\in A$, which is 
written as $\varphi \sim_{\rm asym} \psi $. 
\end{definition}

For AF-algebras, Katsura and Phillips characterized the existence of automorphisms with the Rohlin property. In the following theorem, we show the existence of automorphisms with the Rohlin property for unital simple nuclear TAF-algebras which belongs to the UCT class and with the property (SI). 
Remark that any unital simple nuclear TAF-algebra which belongs to the UCT class is $\Zj$-absorbing (Toms and Winter \cite{TW}). Thus the following theorem applies to these $C^*$-algebras with property (SI).
      
\begin{theorem}\label{thm:1}
Let $A$ be a unital TAF-algebra which satisfies $A \otimes \Zj \cong A$. Suppose that  $*$-homomorphism $\Phi : A \rightarrow A \otimes \Zj $ which is defined by $\Phi(a)=a\otimes 1_{\Zj}$ is asymptotically unitarily equivalent to the isomorphism between $A$ and $A \otimes \Zj$ and $A$ has the property (SI).
Then for any automorphism $\alpha$ of $A$ there exists an automorphism $\widetilde{\alpha} $ of $A$ such that $\widetilde{\alpha} $ has the Rohlin property and $\alpha \sim_{\rm asym}\widetilde{\alpha} $.
\end{theorem}
\bprf
 Let $\alpha$ be an automorphism of $A$ and let $\Psi$ be an isomorphism between $A$ and $A\otimes \Zj $ such that $\Psi \sim_{\rm asym}  \Phi $. Then it follows that 
\begin{eqnarray}
\alpha&=& \Psi^{-1} \circ \Psi \circ \alpha  \sim_{\rm asym}  \Psi^{-1} \circ \Phi \circ \alpha \nonumber \\
&=& \ \Psi^{-1} \circ (\alpha \otimes \sigma) \circ \Phi \sim_{\rm asym} 
\Psi^{-1}  \circ (\alpha \otimes \sigma ) \circ \Psi. \nonumber
\end{eqnarray}
Then by Lemma \ref{lem:p2}, we may assume that
$\alpha $ has a $(p_n) \in P(A_{\infty})$ satisfying the condition in Lemma \ref{lem:p3} for any $k \in \N $.
Note that, for any $u\in U(A) $, $\Ad u \circ \alpha $ equals $\alpha$ on $A_{\infty}$. 
Let $F_n$ be an increasing sequence of finite subsets of $A^1$ such that $\bigcup F_n$ is dense in $A^1$ and equip $I= \{(n,k)\in \N \times \N; k\leq n \}$ with the lexicographic order. We will identify $(n,0)$ with $(n-1,n-1)$ below.

Set $\alpha_{1,1}= \alpha $. We shall construct $\alpha_{(n,k)} \in \Aut (A)$ for
each $(n,k) \in I$ such that $\alpha_{(n,k)} \sim_{\rm asym} \alpha_{(n,k -1)}$ and that $\alpha_{(n,k)}$ has Rohlin towers consisting of $2k+1$ projections with a centrality condition as described in Lemma \ref{lem:p3}. More precisely, by applying Lemma \ref{lem:p3} to $\alpha_{(n,k -1)}$ at each step, we inductively construct the following three objects: a sequence of partitions of unity $\{e_{n,k ,j}^{(i)}\in P(A) ;j=0,1,...,k-1+i,\ i=0,1 \}$, a sequence of unitaries $w_{n,k }$ in $U(C([0,1])\otimes A) $, and  a sequence of finite subsets $G_{n,k}$ of $A^1$ for $(n,k) \in I$ satisfying the following conditions: 
\begin{description}
\item[(I)] $\sum_{i=0,1}\sum_{j=0}^{k-1+i} e_{n,k,j}^{(i)} = 1_A$,

$\|[e_{n,k,j}^{(i)},x]\| <(n+k)^{-1},\quad x \in G_{n,k-1},\ i=0,1,\ j=0,1,...,k-2+i$,

$\displaystyle \|\Ad w_{n,k}(1) \circ \alpha_{n,k-1} (e_{n,k,j}^{(i)}) - e_{n,k,j+1}^{(i)}\|<(n+k)^{-1},\quad \ i,\ j  $,
\item[(I\hspace{-.1em}I) ] $\|[w_{n,k}(t),x] \|<2^{-(n+k)} ,\quad x\in G_{n,k-1},\ i,\ j,\ t\in [0,1],\\ w_{n,k}(0)=1 $,
\item[(I\hspace{-.1em}I\hspace{-.1em}I) ] $\displaystyle G_{n,k} \supset F_n \cup G_{n,k-1} \cup  \alpha_{n,k}(G_{n,k-1})\cup \{e_{n,k,j}^{(i)}\}_{i,j} \cup \{\alpha_{n,k}(e_{n,k,j}^{(i)}) \}_{i,j}$,
\end{description} 
where $\alpha_{n,k} = \Ad w_{n,k}(1) \circ \alpha_{n,k-1}$.
From (I\hspace{-.1em}I) and (I\hspace{-.1em}I\hspace{-.1em}I), we define an automorphism $\widetilde{\alpha}$ of $A$ by
\[ \widetilde{\alpha} (x) = \lim_{(n,k)\rightarrow \infty} \alpha_{n,k} (x). \]
Indeed, for $x\in G_{n,k-1}$ we have $\alpha_{n,k}(x)\approx_{2^{-(n+k)}}\Ad w_{n,k+1}(1)\circ \alpha_{n,k}(x)=\alpha_{n,k+1}(x) \approx_{2^{-(n+k+1)}} \Ad w_{n,k+2}(1) \circ \alpha_{n,k+1}(x)=\alpha_{n,k+2}(x) \approx_{2^{-(n+k+2)}} \cdots$, where $x\approx_{\delta} y $ means that $\|x-y\| < \delta$.  Then it follows that $\alpha_{n,k}(x) $ is a Cauchy sequence. Since $\bigcup G_n $ is a dense subset of $A^1$, we conclude that $\widetilde{\alpha}$ is  an endomorphism of $A$. Similarly, for $x \in G_{n,k}$, we have that 
\[\alpha_{n,k}^{-1}(x) \approx _{2^{-(n+k+1)}} \alpha_{n,k}^{-1}\circ \Ad w_{n,k+1}^*(1) (x) = \alpha_{n,k+1}^{-1} (x) .\]
Then we also obtain  an endomorphism $\widetilde{\alpha}^{-1}$ of $A$ such that 
\[ \widetilde{\alpha}^{-1} (x) = \lim_{(n,k)\rightarrow \infty} \alpha_{n,k}^{-1}(x),\quad x \in \bigcup G_{n,k}.\] We can easily prove that $\widetilde{\alpha}\circ \widetilde{\alpha}^{-1}=\id_A= \widetilde{\alpha}^{-1}\circ \widetilde{\alpha} $ and thus $\widetilde{\alpha}^{-1}$ is indeed the inverse of $\widetilde{\alpha}$ and $\widetilde{\alpha}$ is an automorphism.
 
Define a continuous path of unitaries $w(t)\in U(A)$, $t\in [0,\infty)$ by 
\[w(t)=w_{n,k}(nt-n(n-1)-(k-1))w_{n,k-1}(1)\cdots w_{1,1}(1),\] 
\[(n-1)+n^{-1}(k-1)\leq t \leq (n-1)+n^{-1}k.\] 
Since $w_{n,k}(t)$ almost commutes with $G_{n,k-1}$ for any $t\in [0,1]$, we have that for $x\in G_{n,k}$ 
\[ \widetilde{\alpha} (x) = \lim_{t \rightarrow \infty} \Ad w(t)\circ \alpha(x).\] Thus $\widetilde{\alpha} \sim_{\rm asym} \alpha$.
From (I), when we fix $k$, $(e_{n,k,j}^{(i)})_{n\in \N} $ is a central sequence in $P(A)$.  
From (I) and (I\hspace{-.1em}I\hspace{-.1em}I) we have that
\begin{eqnarray}
e_{n,k,j+1}^{(i)} &\approx _{(n+k)^{-1}}& \Ad w_{n,k} (1) \circ \alpha_{n,k-1}(e_{n,k,j}^{(i)}) = \alpha_{n,k}(e_{n,k,j}^{(i)}) \nonumber \\ 
&\approx_{2^{-(n+k+1)}}& \Ad w_{n,k+1} (1) \circ \alpha_{n,k} (e_{n,k,j}^{(i)}) = \alpha_{n,k+1}(e_{n,k,j}^{(i)}) \nonumber \\
&\approx_{2^{-(n+k+2)}}&\Ad w_{n,k+2} (1) \circ \alpha_{n,k+1} (e_{n,k,j}^{(i)})= \alpha_{n,k+2}(e_{n,k,j}^{(i)}) \approx \cdots .\nonumber
\end{eqnarray}
Then $\| \widetilde{\alpha} (e_{n,k,j}^{(i)})-e_{n,k,j+1}^{(i)} \| < 2(n+k)^{-1}$. Thus we have that $\widetilde{\alpha} ((e_{n,k,j}^{(i)})_n) = (e_{n,k,j+1}^{(i)})_n$ for any $i=0,1$ and any $j=0,1,...,k-2+i$, which shows that $\widetilde{\alpha}$ has the Rohlin property. This completes the proof.
\eprf			

\medskip

{\bf Acknowledgment.} 
The author is grateful to his advisor Akitaka Kishimoto for his support and the many valuable discussions and to Takeshi Katsura and N. Christopher Phillips for a valuable preprint and discussions at an early stage of this research.
This research was partially carried out during the Fields Institute Thematic Program on Operator Algebras in the fall of 2007. The author would like to thank that institution for its support.

\end{document}